\documentclass[12pt]{article}
\usepackage{booktabs}
\usepackage{threeparttable}
\usepackage{float}

\RequirePackage[colorlinks,citecolor=blue,urlcolor=blue]{hyperref}

\def\spacingset#1{\renewcommand{\baselinestretch}%
	{#1}\small\normalsize} \spacingset{1}

\usepackage{amsthm,amsmath,enumerate,amsbsy,amsfonts,amssymb,mathabx,amscd,graphicx,algorithm,  array, booktabs,indentfirst,setspace,subcaption,caption,booktabs,float,enumitem,enumerate}
\usepackage{multirow}
\usepackage{natbib}
\usepackage{geometry}
\usepackage{authblk, caption, subcaption}
\usepackage{mathrsfs}
\usepackage{epsfig}

\newtheorem{theorem}{Theorem}

\newtheorem{remark}{Remark}
\newtheorem{condition}{Condition}

\newcommand{\blambda}{\boldsymbol \lambda}

\newcommand{\by}{{\mathbf y}}

\newcommand{\bW}{{\bf W}}

\newcommand{\bI}{{\bf I}}

\newcommand{\bP}{{\bf P}}
\newcommand{\bS}{{\bf S}}

\newcommand{\bX}{{\bf X}}

\newcommand{\bZ}{{\bf Z}}

\newcommand{\bU}{{\bf U}}
\newcommand{\bV}{{\bf V}}

\newcommand{\bG}{{\bf G}}

\newcommand{\bee}{\boldsymbol{e}}

\newcommand{\var}{\text{Var}}

\newcommand{\diag}{\mbox{diag}}

\newcommand{\cq}{\textnormal{CQ}}
\newcommand{\TB}{\textnormal{TB}}
\newcommand{\DB}{\textnormal{DB}}
\newcommand{\TR}{\textnormal{TR}}
\newcommand{\DR}{\textnormal{DR}}
\newcommand{\SRD}{\textnormal{S}}
\newcommand{\FRD}{\textnormal{F}}

\def\T{{ \mathrm{\scriptscriptstyle T} }}

\makeatletter
\newcommand*{\rom}[1]{\expandafter\@slowromancap\romannumeral #1@}
\makeatother

\makeatletter
\DeclareRobustCommand\widecheck[1]{{\mathpalette\@widecheck{#1}}}
\def\@widecheck#1#2{%
	\setbox\z@\hbox{\m@th$#1#2$}%
	\setbox\tw@\hbox{\m@th$#1%
		\widehat{%
			\vrule\@width\z@\@height\ht\z@
			\vrule\@height\z@\@width\wd\z@}$}%
	\dp\tw@-\ht\z@
	\@tempdima\ht\z@ \advance\@tempdima2\ht\tw@ \divide\@tempdima\thr@@
	\setbox\tw@\hbox{%
		\raise\@tempdima\hbox{\scalebox{1}[-1]{\lower\@tempdima\box
				\tw@}}}%
	{\ooalign{\box\tw@ \cr \box\z@}}}
\makeatother
\RequirePackage[colorlinks,citecolor=blue,urlcolor=blue]{hyperref}

\usepackage{booktabs}
\usepackage{threeparttable}
\usepackage{float}

\usepackage[section]{placeins}

\RequirePackage[colorlinks,citecolor=blue,urlcolor=blue]{hyperref}

\def\spacingset#1{\renewcommand{\baselinestretch}%
	{#1}\small\normalsize} \spacingset{1}

\usepackage{amsthm,amsmath,enumerate,amsbsy,amsfonts,amssymb,mathabx,amscd,graphicx,algorithm,  array, booktabs,indentfirst,setspace,subcaption,caption,booktabs,float,enumitem,enumerate}
\usepackage{multirow}
\usepackage{natbib}
\usepackage{geometry}
\usepackage{authblk, caption, subcaption}
\usepackage{mathrsfs}
\usepackage{epsfig}

\def\T{{ \mathrm{\scriptscriptstyle T} }}



\makeatletter
\DeclareRobustCommand\widecheck[1]{{\mathpalette\@widecheck{#1}}}
\def\@widecheck#1#2{%
	\setbox\z@\hbox{\m@th$#1#2$}%
	\setbox\tw@\hbox{\m@th$#1%
		\widehat{%
			\vrule\@width\z@\@height\ht\z@
			\vrule\@height\z@\@width\wd\z@}$}%
	\dp\tw@-\ht\z@
	\@tempdima\ht\z@ \advance\@tempdima2\ht\tw@ \divide\@tempdima\thr@@
	\setbox\tw@\hbox{%
		\raise\@tempdima\hbox{\scalebox{1}[-1]{\lower\@tempdima\box
				\tw@}}}%
	{\ooalign{\box\tw@ \cr \box\z@}}}
\makeatother
\RequirePackage[colorlinks,citecolor=blue,urlcolor=blue]{hyperref}

\usepackage{bookmark}

\setcitestyle{aysep={,},yysep={;}}

\allowdisplaybreaks
\numberwithin{remark}{section}

\makeatletter
\newcommand{\figcaption}{\def\@captype{figure}\caption}
\newcommand{\tabcaption}{\def\@captype{table}\caption}
\makeatother

\makeatletter
\DeclareRobustCommand\widecheck[1]{{\mathpalette\@widecheck{#1}}}
\def\@widecheck#1#2{%
	\setbox\z@\hbox{\m@th$#1#2$}%
	\setbox\tw@\hbox{\m@th$#1%
		\widehat{%
			\vrule\@width\z@\@height\ht\z@
			\vrule\@height\z@\@width\wd\z@}$}%
	\dp\tw@-\ht\z@
	\@tempdima\ht\z@ \advance\@tempdima2\ht\tw@ \divide\@tempdima\thr@@
	\setbox\tw@\hbox{%
		\raise\@tempdima\hbox{\scalebox{1}[-1]{\lower\@tempdima\box
				\tw@}}}%
	{\ooalign{\box\tw@ \cr \box\z@}}}
\makeatother
\RequirePackage[colorlinks,citecolor=blue,urlcolor=blue]{hyperref}

\usepackage{geometry}
\newcommand{\blind}{1}

\def\spacingset#1{\renewcommand{\baselinestretch}%
{#1}\small\normalsize} \spacingset{1}

\begin{document}
	\if1\blind
	{
		\title{\bf \Large On Robust Empirical Likelihood for Nonparametric Regression with Application to Regression Discontinuity Designs }
		\author[1]{Qin Fang}
        \author[2]{Shaojun Guo}
           \author[3]{Yang Hong}
   \author[4]{Xinghao Qiao} 
        \affil[1]{University of Sydney Business School, Sydney, Australia}
		\affil[2]{Institute of Statistics and Big Data, Renmin University of China, China}
        \affil[3]{School of Mathematics, Renmin University of China, China}
		\affil[4]{Faculty of Business and Economics, The University of Hong Kong, Hong Kong}
  \setcounter{Maxaffil}{0}
		
		\renewcommand\Affilfont{\itshape\small}
		\date{\vspace{-3ex}}
		\maketitle
	} \fi
	\if0\blind
	{
		\bigskip
		\bigskip
		\bigskip
		\begin{center}
			{\Large\bf On Robust Empirical Likelihood for Nonparametric Regression with Application to Regression Discontinuity Designs }
		\end{center}
		\medskip
	} \fi

\spacingset{1.25}

\begin{abstract}

Empirical likelihood serves as a powerful tool for constructing confidence intervals in nonparametric regression and regression discontinuity designs (RDD). The original empirical likelihood framework can be naturally extended to these settings using local linear smoothers, with Wilks' theorem holding only when an undersmoothed bandwidth is selected. However, the generalization of bias-corrected versions of empirical likelihood under more realistic conditions is non-trivial and has remained an open challenge in the literature. This paper provides a satisfactory solution by proposing a novel approach, referred to as robust empirical likelihood, designed for nonparametric regression and RDD. The core idea is to construct robust weights which simultaneously achieve bias correction and account for the additional variability introduced by the estimated bias, thereby enabling valid confidence interval construction without extra estimation steps involved. We demonstrate that the Wilks' phenomenon still holds under weaker conditions in nonparametric regression, sharp and fuzzy RDD settings. Extensive simulation studies confirm the effectiveness of our proposed approach, showing superior performance over existing methods in terms of coverage probabilities and interval lengths. Moreover, the proposed procedure exhibits robustness to bandwidth selection, making it a flexible and reliable tool for empirical analyses. The practical usefulness is further illustrated through applications to two real datasets.
\end{abstract}

\bigskip
\noindent%


\newpage
\spacingset{1.45}

\setlength{\abovedisplayskip}{0.2\baselineskip}
\setlength{\belowdisplayskip}{0.2\baselineskip}
\setlength{\abovedisplayshortskip}{0.2\baselineskip}
\setlength{\belowdisplayshortskip}{0.2\baselineskip}

\section{Introduction}
\label{sec:intro}


{\color{black}

Nonparametric regression has emerged as a key tool in statistics and econometrics for modeling complex relationships between dependent and independent variables without imposing rigid functional-form restrictions. 
Empirical likelihood (EL), introduced by \cite{owen1988empirical,Owen1990}, provides likelihood-based inference without requiring a fully specified parametric model. Among its appealing features, EL produces naturally shaped confidence sets that respect parameter constraints, admit Bartlett correction, and enjoy Wilks-type, distribution-free limiting distributions \citep{Owen2001book}.
When combined with local polynomial smoothing \citep{fan1996local},  however, EL inherits a central difficulty of nonparametric inference.
The bandwidths chosen for estimation are often too large for confidence interval construction, leaving a non-negligible smoothing bias that can distort the Wilks approximation for the EL ratio.
 In this paper, we aim to develop 
 a robust and  easy-to-implement
 bias-corrected EL inference procedure for nonparametric regression, with an application to regression discontinuity designs.
 Our method automatically restores the appeal of EL as a likelihood-based, asymptotically distribution-free inference tool and delivers confidence intervals with close-to-nominal empirical coverage across a broader range of bandwidth choices, thereby providing a clean and satisfactory solution to a long-standing challenge in the EL literature.

Existing EL-based inferences via local linear smoothers primarily fall into two categories.  One type of approach is based on undersmoothing, where the smoothing bandwidth $h$ is chosen to be smaller than the optimal value
such that the bias becomes asymptotically negligible relative to the variance; see, e.g., \cite{chenQin2000}, \cite{otsu2012empirical} and \cite{Otsu2015}. While producing valid coverage asymptotically, undersmoothing often results in increased variance and overly wide confidence intervals in finite samples.

The other approach relies on bias correction techniques, which involve a pilot bandwidth
$b$ to estimate and remove the asymptotic bias. In the Wald-type literature, the robust bias-corrected inference framework of \cite{Calonico2014} estimates and removes the leading bias while calibrating inference to account for the additional variability introduced by bias estimation. Subsequent work further developed the coverage properties, bandwidth selection, and coverage-error-optimal implementation of this approach; see \cite{calonico2018jasa,calonico2020ej,calonico2022coverage}. However, 
the analogous problem for EL is more delicate, as bias correction changes the estimating equations, and the variability introduced by the estimated bias can alter the Wilks-type calibration of the EL ratio.  
 To examine this issue, we investigate, in Section~\ref{sec:pitfalls}, two conventional bias-corrected EL methods in nonparametric regression, where the bias estimators are built upon Taylor expansion and direct differencing, respectively.
Under mild regularity conditions, we show that the asymptotic distributions of the resulting EL ratios depend explicitly on the ratio $h/b$. The standard chi-squared limit holds only
when $h/b \to 0$ as $n\to \infty$; otherwise, the Wilks' phenomenon fails. 
However, such a ratio condition is often violated in practice, as noted by \cite{Calonico2014}. Consequently, confidence intervals formed from these conventional bias-corrected EL ratios using standard chi-squared cutoffs often undercover, as demonstrated in Section \ref{sec:sim} through extensive simulations. 
This then leaves open a key challenge in the EL literature, namely incorporating bias correction into EL while retaining Wilks-type calibration under the practically relevant regime \(h/b\to\kappa\in(0,1]\).

To address this problem, we propose a novel bias-corrected EL framework, called robust EL.
Unlike existing methods,
which typically handle first-order impact by adjusting EL ratios, 
our proposal adopts a fundamentally different strategy
by developing new sets of weights, termed ``robust weights", which simultaneously achieve bias correction and account for the additional variability introduced by the bias estimation in the EL formulation.
Building on the new weighting schemes in Section~\ref{sec:3.1}, we propose two types of robust EL ratio functions and demonstrate that, under standard regularity conditions in nonparametric regression, both robust EL ratios converge to a standard chi-squared distribution within the refined asymptotic framework of \cite{Calonico2014}, where the ratio 
\begin{equation}
\label{eq.ratio}
    \frac{h}{b} \to \kappa \in [0, 1], \text{  as  } n \to \infty.
\end{equation} 
In fact, the term ``robust" in our proposal is inspired by the robust bias correction framework of \cite{Calonico2014}. While our method can also be described as robust weighted EL or robust bias-corrected EL, we refer to it simply as robust EL for brevity.

Regression discontinuity designs (RDD) provide a natural and important application of the proposed framework. As a leading quasi-experimental method in causal inference, RDD identifies treatment effects from discontinuities in treatment assignment at a cutoff, where nonparametric local polynomial methods are routinely used for inference; see \cite{cattaneo2022regression} for a recent review. EL-based inference for sharp and fuzzy RDD settings was first studied by \cite{Otsu2015}, relying on undersmoothing to obtain asymptotically valid confidence intervals. More recently, \cite{ma2020empirical} proposed a scale-adjusted EL ratio to restore the pivotal property asymptotically. Although sometimes effective, the scale factor depends on estimating variance and covariance of the local linear estimates and the estimated bias, making it cumbersome to implement and diminishing the distribution-free appeal of EL.

We thus further extend the proposed robust EL approaches to 
accommodate both sharp and fuzzy RDD settings, establishing the corresponding versions of Wilks' Theorem in this context. Numerical studies in Section \ref{sec:sim} confirm that our methods consistently outperform the original EL based on undersmoothing, conventional bias-corrected EL procedures, and the normal approximation-based approach of \cite{Calonico2014}, particularly in terms of empirical size and coverage probabilities.

Our paper makes three main contributions. Methodologically,  we pioneer a bias-corrected EL inference framework for nonparametric regression based on fully data-adaptive robust weights.
The proposed weights jointly correct bias and address its variability, enabling a single-pass computation of robust confidence intervals. Theoretically, our method preserves the distribution-free nature of EL and its Wilks-type calibration, allowing confidence intervals to be constructed using standard chi-squared critical values under the more general ratio condition in \eqref{eq.ratio}. This avoids undersmoothing and provides robustness over a broader range of choices for \(h\) and \(b\).

On the application side, 
we show that the proposed robust EL procedure and the associated theory can be seamlessly extended to both sharp and fuzzy RDD settings by adapting the robust weights and EL constraints accordingly.  Under some additional conditions, the pivotal property continues to hold asymptotically. Beyond first-order validity, we also explore the second-order asymptotic properties of the proposed EL methods and identify a potential route toward attaining a sharp coverage-error rate.
More broadly, our robust EL framework is potentially applicable to quantile RDD and other scenarios; we discuss these extensions further in Section~\ref{sec:diss}. 


The rest of the paper is organized as follows. In Section~\ref{sec:pitfalls}, we discuss the limitations of original and conventional bias-corrected  EL methods in nonparametric regression. Section~\ref{sec:methodology} presents the construction of robust weights and details our proposed robust EL approaches, along with their asymptotic distributions. In Section~\ref{sec:RDD}, we apply our proposal to both sharp and fuzzy RDD and establish the corresponding theoretical validity. Sections~\ref{sec:sim} and \ref{sec:real} illustrate the superiority of our robust EL approaches over the original and conventional bias-corrected EL approaches through extensive simulation studies and the analysis of two real datasets, respectively.
Section~\ref{sec:diss} concludes the paper and discusses several extensions. All technical proofs are relegated to the appendix.
}

\section{Pitfalls of conventional EL methods}
\label{sec:pitfalls}

To facilitate the construction of our robust confidence intervals for nonparametric regression, we  start with  an overview of local polynomial estimation and the original EL-based inference procedure of \cite{chenQin2000} in Section~\ref{sec:2.1}, and present two conventional attempts at bias-corrected EL-based inference with asymptotic analysis in Section~\ref{sec:2.2}.

\subsection{Methodological background}
\label{sec:2.1}
Suppose we observe $\{(X_i, Y_i)\}_{i = 1}^n$, a set of $n$ independent pairs from the nonparametric regression model, 
\begin{equation} \label{eq1}
Y_i = m(X_i) + \varepsilon_i,\quad i=1, \dots,n,
\end{equation}
where $X_i$ follows a density function $f(\cdot)$ with bounded support $[0,1]$, and $\{\varepsilon_i\}$ is a sequence of independent errors satisfying $E(\varepsilon_i | X_i = x) = 0$ and $\sigma^2(x) = E(\varepsilon_i^2 | X_i = x)<\infty$. 
In what follows, denote $K_h(\cdot)=K(\cdot/h)/h$ for a univariate kernel $K$ with bandwidth $h>0.$
We define the design matrix  of order $p$ as
\begin{equation*}
\bX_{p,h} = \left(
\begin{array}{cccc}
1 & \frac{X_1 - x}{h} & \cdots & \left(\frac{X_1 - x}{h}\right)^p \\
\vdots & \vdots & & \vdots \\
1 & \frac{X_n - x}{h} & \cdots & \left(\frac{X_n - x}{h}\right)^p \\
\end{array}
\right)~~{\in \mathbb R}^{n \times (p+1)}.
\end{equation*}
Let $\bW_h = \text{diag}\{K_h(X_i - x)\}$ be the $n \times n$ diagonal matrix of kernel weights. 
We write $S_{j,h}(x) = n^{-1} \sum_{i=1}^{n} K_h(X_i - x) (X_i - x)^j / h^j$, and put $\bS_{p,h} = n^{-1} \bX_{p,h}^\T \bW_h \bX_{p,h}$, whose $(j,k)$-th entry is given by $S_{j+k-2, h}$ for $j,k=1, \dots, p+1$.

The local polynomial estimator of order $p$ for the $\ell$-th derivative $m^{(\ell)}(x)$ with bandwidth $h$,  denoted as $\widehat{m}^{(\ell)}_{p,h}(x)$, at each $x \in [0,1]$, takes the form of
\begin{equation*}
\begin{aligned}
\widehat{m}^{(\ell)}_{p,h}(x) = \ell! (nh^\ell)^{-1} \bee_{\ell+1}^\T \bS_{p,h}^{-1} \bX_{p,h}^\T \bW_h \by 
= \frac{\ell!}{nh^\ell} \sum_{i=1}^n W_{i,\ell,p,h}(x) Y_i, \quad \ell = 0, \ldots, p,
\end{aligned}
\end{equation*}
where $\bee_{\ell+1} = (0, \ldots, 0, 1, 0, \ldots, 0)^\T$ is a unit vector with the $1$ in the $(\ell+1)$-th position, $\by = (Y_1, \ldots, Y_n)^\T$, and $W_{i,\ell,p,h}(x) = \bee_{\ell+1}^\T \bS_{p,h}^{-1} \big\{1, (X_i - x)/h, \ldots, (X_i - x)^p/h^p\big\}^\T K_h(X_i - x)$ for $i=1, \dots, n.$
We define the local linear weights
as $W_{i,h}(x) = K_h(X_i - x) \big\{S_{2,h}(x) - S_{1,h}(x) (X_i - x)/h\big\}$.
In the special case of $\ell = 0$, the local linear estimator of $m(x)$ with bandwidth $h$ then simplifies to
\[
\widehat{m}_{1,h}(x) = \frac{1}{n} \sum_{i=1}^n W_{i,0,1,h}(x) Y_i = \frac{\sum_{i=1}^n W_{i,h}(x) Y_i}{\sum_{i=1}^n W_{i,h}(x)}.
\]
Note that  $\widehat{m}_{1,h}(x)$ satisfies the weighted moment equation 
$
\sum_{i=1}^n W_{i,h}(x) \left (Y_i - \theta\right) = 0.
$
Inspired by this, \cite{chenQin2000} (referred to as CQ) introduced the 
 original empirical log-likelihood ratio for $m(x)$, given by
\begin{equation*}
\begin{aligned}
l_{\cq}(\theta) = -2 \max \bigg\{ \sum_{i=1}^n \log(np_i) \,\bigg|\, p_i \ge 0, \,\sum_{i=1}^n p_i = 1, \,\sum_{i=1}^n p_i W_{i,h}(x) (Y_i - \theta) = 0 \bigg\}.
\end{aligned}
\end{equation*}
However, this EL-based inference targets $m(x) + \text{bias}(x),$ rather than $m(x)$ itself. Under standard regularity conditions,
\cite{chenQin2000} demonstrated that $l_{\cq}(m(x))$ asymptotically follows a standard chi-square distribution with one degree of freedom, provided that $nh^5 \to 0$ and $m^{(2)}(x) \neq 0$. Hence, undersmoothing is required to reduce the impact of bias in this formulation.


\subsection{Two conventional bias-corrected EL methods}
\label{sec:2.2}

We modify the original EL-based inference procedure for $m(x)$ by directly incorporating two commonly-used bias correction strategies. 
The first approach is motivated by 
 the Taylor expansion of $m(X_i)$ around $x$.  
 Let $r(X_i) = m(X_i) - m(x) - m^{(1)}(x)(X_i - x)$ denote the approximation error at $X_i$. Applying a second-order expansion, we can approximate this error as $r(X_i) \approx m^{(2)}(x)(X_i - x)^2/2.$
Since the bias
in  $\widehat{m}_{1,h}(x)$ arises from $r(X_i)$ near $x$, we consider a bias estimator as
\begin{equation} \label{r_1b}
    \widehat{r}_{1,b}(X_i) = \frac{1}{2}\widehat{m}_{2,b}^{(2)}(x)(X_i - x)^2,
\end{equation}
where $
\widehat{m}_{2,b}^{(2)}(x) = 2\sum_{i=1}^n W_{i,2,2,b}(x) Y_i /(nb^2)
$ denotes the local quadratic estimator of $m^{(2)}(x)$ with
 a pilot bandwidth $b$.  
See \cite{fan1998local} for a detailed discussion on bias correction in general nonparametric settings and \cite{Calonico2014} for its applications in RDD. 
 The corresponding {\bf T}aylor-expansion-based {\bf B}ias-corrected (TB) empirical log-likelihood ratio for $m(x)$ is  then defined as
\begin{equation}
\label{bc-1}
\begin{aligned}
l_{\TB}(\theta) = - 2 \max \bigg\{\sum_{i = 1}^n \log(np_i) \,\bigg|\, & p_i \ge 0, \, \sum_{i=1}^n p_i = 1,\\
& \sum_{i=1}^n p_iW_{i,h}(x) \big\{Y_i - \widehat{r}_{1,b}(X_i)- \theta\big\} = 0 \bigg\}.
\end{aligned}
\end{equation}

The second bias-corrected approach is inspired by the direct difference  $m(X_i) - m(x)$.  Rather than relying on a Taylor expansion, we approximate the bias using the local difference estimation, i.e.,
$$
\widehat{r}_{2,b}(X_i) = \widehat{m}_{2,b}(X_i) - \widehat{m}_{2,b}(x),
$$ 
 where $\widehat{m}_{2,b}(x)$ represents the local quadratic estimator of $m(x)$ with a pilot bandwidth $b>0$.
Substituting this into the EL formulation, we obtain 
the  {\bf D}ifference-based {\bf B}ias-corrected (DB) empirical log-likelihood ratio, defined as
\begin{equation}
\label{bc-2}
\begin{aligned}
l_{\DB}(\theta) = - 2 \max \bigg\{\sum_{i = 1}^n \log(np_i) \,\bigg|\, & p_i \ge 0, \,\sum_{i=1}^n p_i = 1,\\
& \sum_{i=1}^n p_i W_{i,h}(x) \big\{Y_i - \widehat{r}_{2,b}(X_i) - \theta\big\} = 0 \bigg\}.
\end{aligned}
\end{equation}
In particular, when the estimated bias $\widehat r_{2,b}(X_i)$  is replaced by  its counterpart based on local linear estimators with bandwidth $h=b$, i.e.,
$
\widehat{m}_{1,h}(X_i) - \widehat{m}_{1,h}(x)
$, 
this reduces to the special case of the so-called residual-adjusted EL method proposed by \cite{xue2007empirical, xue2007empirical2} for longitudinal data. Despite the different forms of bias correction, our established Theorem~\ref{th1} below remains applicable to their case.

  

To investigate the theoretical properties, we impose some regularity conditions.

\begin{condition}\label{cond1}
 $E(\varepsilon_i^{4}) < \infty$. 
	\end{condition}
    
\begin{condition}\label{cond2}
{\rm (i)} Both $f(x)$, the density function of $X$, and $\sigma^2(x) = E(\varepsilon_i^2 | X_i = x)$ are continuous and bounded away from zero; {\rm (ii)} The function $m(\cdot)$ is thrice continuously differentiable in a neighborhood of $x$.
\end{condition}

\begin{condition}\label{cond3}
The kernel function $K(\cdot)$ is a symmetric and bounded density function with support $[-1,1]$.
\end{condition}

\begin{condition}\label{cond4}
As $n \to \infty,$  $h \to 0,$  $b^2\log n \to 0,$ $nh^5b^2 \to 0$, $nh^3b^4 \to 0$, $nh \to \infty$, and $n^{\epsilon -1} b \to \infty$ with $1 < \epsilon < 3/2$.
\end{condition}

Conditions~\ref{cond1}--\ref{cond4} are standard in the nonparametric regression literature, see, e.g., \cite{fan1999,chenQin2000,Calonico2014}. The smoothness requirement on $m(\cdot)$ in Condition~\ref{cond2} can be  relaxed to being twice continuously differentiable, in which case the leading bias term reduces to $o(h^2)$ instead of $O(h^3)$. For simplicity in our proofs, we choose to adopt the stronger condition here. The requirements $nh^5b^2 \to 0$ and $nh^3b^4 \to 0$ in Condition~\ref{cond4} ensure that both bias-correction terms are asymptotically unbiased. 

Before presenting the theoretical results, we introduce a refined asymptotic framework that characterizes a more flexible and practically reasonable relationship between the smoothing bandwidth $h$ and the pilot bandwidth $b$.  For local polynomial estimation, balancing bias and variance typically results in the optimal bandwidth selections $h \asymp n^{-1/5}$ for $\widehat m_{1,h}(x)$ and $b \asymp n^{-1/7}$ for $\widehat{m}_{2,b}^{(2)}(x)$. Given these rates, it is generally accepted in the literature to assume that the ratio $h/b\to 0$.  While this assumption is theoretically convenient, determining how small $h/b$ should be in practice remains unclear. For instance, when $h = 0.1$ and $b = 0.3$, the ratio $1/3$ may not be sufficiently small to effectively reduce the first-order impact of the bias correction. 
Thus, we instead follow 
\cite{Calonico2014} and impose the refined asymptotic framework in (\ref{eq.ratio}).
We now proceed to examine the asymptotic distributions of these two conventional bias-corrected EL ratio functions.


\begin{theorem}
\label{th1}
Suppose that Conditions {\rm \ref{cond1}--\ref{cond4}} hold, $m(x)$ is the true parameter, and $x$ is an interior point.
\begin{itemize}
\item[(a)] If $h/b \to 0$ as $n \to \infty,$ then, 
\begin{equation*}
l_{\TB}\left\{m(x)\right\} \stackrel{\mathcal{D}}{\rightarrow} \chi_{1}^{2},
~~\mbox{and}~~l_{\DB}\left\{m(x)\right\} \stackrel{\mathcal{D}}{\rightarrow} \chi_{1}^{2},
\end{equation*}
as $n \to \infty,$ where $\stackrel{\mathcal{D}}{\rightarrow}$ means the convergence in distribution;
\item[(b)] If $h/b \to \kappa \in (0,1]$  as $n \to \infty,$ then, 
\begin{equation*}
l_{\TB}\left\{m(x)\right\} \stackrel{\mathcal{D}}{\rightarrow} \gamma_{1,\kappa}\chi_{1}^{2},
~~\mbox{and}~~l_{\DB} \left\{m(x)\right\} \stackrel{\mathcal{D}}{\rightarrow} \gamma_{2,\kappa}\chi_{1}^{2},
\end{equation*}
as $n \to \infty,$ where the constants $\gamma_{1,\kappa}$ and $\gamma_{2,\kappa}$ are specified in Section~{\rm \ref{appendix-1}} in the appendix  with $\gamma_{1,\kappa} \neq 1$ and $\gamma_{2,\kappa} \neq 1$.
\end{itemize}
\end{theorem}

\begin{remark}
\label{rm.1}
Theorem~\ref{th1} highlights the critical role of the bandwidth ratio $h/b$ in shaping the asymptotic properties of the conventional bias-corrected EL ratios.
In particular, when $\kappa \neq 0$, these EL ratios asymptotically deviate from the standard chi-squared distribution. Instead, their limiting distributions involve nuisance parameters, i.e., $\gamma_{1,\kappa} \neq 1$ and $\gamma_{2,\kappa} \neq 1$, which complicates their use for inference. If this nonstandard behavior is ignored,  hypothesis tests based on  the $\chi_1^2$ distribution may substantially over-reject true null hypotheses, leading to inflated Type I error rates. Likewise, the resulting confidence intervals may undercover, leading to invalid inference conclusions. One possible alternative is to estimate \(\gamma_{1,\kappa}\) or \(\gamma_{2,\kappa}\) and use an appropriately scaled chi-squared distribution. However, such an adjustment not only increases the implementation complexity but, more importantly, compromises the distribution-free nature of empirical likelihood, making it a less desirable solution in both theory and practice.
\end{remark}

\begin{remark}
 \textcolor{black}{\citet{Hjort2009} developed a general plug-in EL theory for inference based on estimating equations with estimated nuisance parameters. A key message of their work is that the usual Wilks phenomenon for EL is not automatically preserved after plug-in estimation, as shown in Theorem ~\ref{th1}. Instead, the limiting distribution of the EL ratio can differ from the standard chi-squared distribution, with the departure reflecting the additional variability introduced by nuisance estimation. Theorem~\ref{th1} provides a local nonparametric analogue of this phenomenon for the conventional bias-corrected EL. In particular, when \(h/b\to\kappa\in(0,1]\), Theorem~\ref{th1} explicitly identifies the first-order limiting distribution as a scaled chi-squared distribution. Thus, as in the plug-in theory of \citet{Hjort2009}, the additional estimation step changes the asymptotic calibration of the EL ratio.}
\end{remark}

\begin{remark}
While Theorem~\ref{th1}  focuses on the case where $x$
 lies in the interior of the support, a similar result holds when $x$ is near the boundary.  The discussion in Remark~\ref{rm.1} then naturally extends to the RDD analysis, where inference at the cutoff is of primary interest. We will explore the RDD settings further in Section~\ref{sec:RDD}.
\end{remark}

The failure of Wilks' theorem when \( h/b \to \kappa \in (0,1] \) can be explained as follows. Define 
$
Z_i(x) = W_{i,h}(x)\left\{Y_i - \widehat{r}_{1,b}(X_i) - m(x)\right\},
$
and consider the following quantities:
\[
U_1(x) = \frac{1}{n}\sum_{i=1}^n Z_i(x),~~U_2(x) = \frac{1}{n}\sum_{i=1}^n \big\{Z_i(x)\big\}^2.
\]
As shown in Section~\ref{appendix-1-1} of the appendix, we can establish that
\[
l_{\TB}\{m(x)\} = \big\{U_2(x)\big\}^{-1}\big\{U_1(x)\big\}^2 \big\{1 + o_P(1)\big\}.
\]
The key insight is that while deriving the limiting distribution of \( U_1(x) \) requires accounting for the additional variability introduced by the bias estimator \( \widehat{r}_{1,b}(X_i) \), the convergence in probability of \( U_2(x) \) depends only on the consistency of the local quadratic estimator \( \widehat{m}_{2,b}^{(2)}(x) \) used in the bias estimation. This indicates that the conventional bias correction formulation incorporates only part of the additional variability,
resulting in a scaled chi-squared limiting distribution. 
A similar limitation applies to \( l_{\DB}\{m(x)\} \). To address this issue, we propose a new strategy in Section~\ref{sec:methodology} by constructing new weights to 
fully capture the variability associated with the bias estimators in a simple and effective way. 

\section{Methodology}
\label{sec:methodology}
In this section, we develop a new EL framework in the nonparametric regression setting, referred to as robust EL. We first propose two sets of robust weights in Section~\ref{sec:3.1} and then construct the corresponding robust EL-based confidence intervals for $m(x)$ with theoretical validity in Section~\ref{sec:3.2}.

\subsection{Construction of robust weights}
\label{sec:3.1}

We begin by revisiting equation (\ref{bc-1}). 
Recall that 
\[
\widehat{m}_{2,b}^{(2)}(x) = \frac{2}{nb^2} \sum_{i=1}^n W_{i,2,2,b}(x) Y_i = \frac{2}{b^2} \sum_{i=1}^n \frac{1}{n}W_{i,2,2,b}(x) (Y_i - \theta),
\]
where the last equality follows from the fact that $\sum_{i=1}^n W_{i,2,2,b}(x) = 0$. 
We first define a `weighted' version of $\widehat{m}_{2,b}^{(2)}(x)$ as
\[
\widehat{m}^{(2)*}_{2,b}(x) = \frac{2}{b^2} \sum_{i=1}^n p_i W_{i,2,2,b}(x) (Y_i - \theta),
\]
where we use weights $\{p_i\}_{i=1}^n$ instead of $n^{-1}$. Substituting $\widehat{m}_{2,b}^{(2)}(x)$ in \eqref{r_1b} with $\widehat{m}^{(2)*}_{2,b}(x)$ and rewriting the moment constraint in (\ref{bc-1}), we obtain that 
\begin{equation*}
\begin{aligned}
&\sum_{i=1}^n p_i W_{i,h}(x) \Big\{Y_i - \frac{1}{2} \widehat{m}^{(2)*}_{2,b}(x) (X_i - x)^2 - \theta\Big\} \\
= &\sum_{i=1}^n p_i \Big\{ W_{i,h}(x) - W_{i,2,2,b}(x) \frac{1}{b^2} \sum_{k=1}^n p_k W_{k,h}(x) (X_k - x)^2 \Big\} (Y_i - \theta) \\
= &\sum_{i=1}^n p_i \Big[ W_{i,h}(x) - W_{i,2,2,b}(x) \frac{1}{nb^2} \sum_{k=1}^n W_{k,h}(x) (X_k - x)^2 \big\{1 + o_P(1)\big\} \Big] (Y_i - \theta),
\end{aligned}
\end{equation*}
where the second equality follows from the approximation
\[
\frac{1}{b^2} \sum_{k=1}^n p_k W_{k,h}(x) (X_k - x)^2 = \frac{1}{nb^2} \sum_{k=1}^n W_{k,h}(x) (X_k - x)^2\big\{1 + o_P(1)\big\}.
\]
Therefore, instead of relying on both original weights $\{W_{i,h}(x)\}_{i=1}^n$ and the bias estimators $\{\widehat{r}_{1,b}(X_i)\}_{i=1}^n$ in (\ref{bc-1}), we introduce the Taylor-expansion-based {\it robust weights} $\{W_{i,h,b}^\star(x)\}_{i=1}^n:$ 
\begin{equation}
\label{rbc-weight}
W_{i,h,b}^\star(x) = W_{i,h}(x) - W_{i,2,2,b}(x) \frac{1}{nb^2} \sum_{k=1}^n W_{k,h}(x) (X_k - x)^2,
\end{equation}
which are fully data-adaptive and better equipped to handle the previously unaddressed uncertainty of $\widehat{r}_{1,b}(X_i)$. 
Some standard calculations yield that
\[
\frac{1}{nb^2} \sum_{k=1}^n W_{k,h}(x) (X_k - x)^2 = O_P(\kappa_n^2)~~\text{with}~~\kappa_n = h/b.
\]
If $\kappa_n \to 0$, then $W_{i,h,b}^\star(x) - W_{i,h}(x) = o_P(1)$, implying that the robust weights asymptotically align with the original ones. However, if $\kappa_n$ approaches a positive constant, then $W_{i,h,b}^\star(x) - W_{i,h}(x) = O_P(1)$, indicating a substantial difference between the robust and original weights. This distinction is crucial in ensuring that the robust bias-corrected EL-based inference procedure, which will be introduced in Section~\ref{sec:3.2},
remains both theoretically valid and practically reliable,  even when $h/b$ does not vanish.
See Theorem~\ref{th2} in Section~\ref{sec:3.2} for further theoretical details.

The same principle can be applied to equation (\ref{bc-2}). Recall that
\[
\widehat{m}_{2,b}(X_i) - \widehat{m}_{2,b}(x) = \frac{1}{n} \sum_{k=1}^n \big\{W_{k,0,2,b}(X_i) - W_{k,0,2,b}(x)\big\} (Y_k - \theta),
\]
where the equality holds as $\sum_{k=1}^n \left\{W_{k,0,2,b}(X_i) - W_{k,0,2,b}(x) \right\} = 0.$ We can write its weighted version as
\[
\widehat{m}_{2,b}^*(X_i) - \widehat{m}_{2,b}^*(x) = \sum_{k=1}^n p_k\big\{W_{k,0,2,b}(X_i) - W_{k,0,2,b}(x)\big\} (Y_k - \theta).
\]
Replacing this weighted estimator in  (\ref{bc-2}) and rewriting the moment constraint  yields that
\begin{equation*}
\begin{aligned}
&\sum_{i=1}^n p_i W_{i,h}(x) \big[Y_i - \left\{\widehat{m}_{2,b}^*(X_i) - \widehat{m}_{2,b}^*(x)\right\} - \theta\big] \\
 = &\sum_{i=1}^n p_i \Big[ W_{i,h}(x) - \sum_{k=1}^n p_k W_{k,h}(x) \big\{W_{i,0,2,b}(X_k) - W_{i,0,2,b}(x)\big\} \Big] (Y_i - \theta) \\
 = &\sum_{i=1}^n p_i \Big[ W_{i,h}(x) - \frac{1}{n} \sum_{k=1}^n W_{k,h}(x) \big\{W_{i,0,2,b}(X_k) - W_{i,0,2,b}(x)\big\} \big\{1 + o_P(1)\big\} \Big] (Y_i - \theta),
\end{aligned}
\end{equation*}
where the second equality is due to
\[
\sum_{k=1}^n p_k W_{k,h}(x) \big\{W_{i,0,2,b}(X_k) - W_{i,0,2,b}(x) \big\}
= \frac{1}{n} \sum_{k=1}^n W_{k,h}(x) \big\{W_{i,0,2,b}(X_k) - W_{i,0,2,b}(x) \big \}\big\{1 + o_P(1)\big\}.
\]
We then define the difference-based {\it robust weights} $\{W_{i,h,b}^{\diamond}(x)\}_{i=1}^n$ as 
\begin{equation}
\label{rbc-weight2}
W_{i,h,b}^{\diamond}(x)  = W_{i,h}(x) - \frac{1}{n}\sum_{k=1}^n W_{k,h}(x) \big \{W_{i,0,2,b}(X_k) - W_{i,0,2,b}(x)\big\},
\end{equation}
which can also automatically adjust for the inherent variability of $\widehat{r}_{2,b}(X_i)$ in (\ref{bc-2}). Similarly, when $\kappa_n \to 0$, the difference between the robust and original weights is negligible, i.e., $W_{i,h,b}^{\diamond}(x) - W_{i,h}(x) = o_P(1)$. When $\kappa_n$ tends to a positive constant $\kappa$, then $W_{i,h,b}^{\diamond}(x) - W_{i,h}(x) = O_P(1)$.
It is worth noting that our numerical experiments demonstrate the superior performance of the robust weights in \eqref{rbc-weight2} over their local linear counterparts, where $W_{i,0,2,b}(X_k)$ and $W_{i,0,2,b}(x)$ are replaced by $W_{i,0,1,b}(X_k)$ and $W_{i,0,1,b}(x)$, respectively. We thus adopt the local quadratic formulation in subsequent analysis.

\begin{remark} \label{rmk.3.1}
{\color{black} The intuition behind the formulations in \eqref{rbc-weight} and
\eqref{rbc-weight2} comes from the fact that the weights $p_i$'s in the EL framework can be interpreted as resampling probabilities. To fully incorporate the variability of the bias estimator into the EL formulation, these resampling probabilities must be explicitly integrated into the additional estimators, namely $\widehat{m}^{(2)}_{2,b}(x)$ in the
Taylor-expansion-based approach, and $\widehat{m}_{2,b}(X_i)$ and
$\widehat{m}_{2,b}(x)$ in the difference-based approach.  From this perspective, our proposed robust EL procedure is conceptually analogous to the double-resampling approach of \cite{guo2021better}, which was  developed to correct variance estimation in the presence of plug-in nuisance estimators.
} 
\end{remark}
\begin{remark} \label{rmk.3.2}
    {\color{black}It is worth noting that under the local-linear and
local-quadratic specification considered here, our
Taylor-expansion-based robust weights in \eqref{rbc-weight} coincide,
up to a multiplicative factor, with the finite-sample weights
$\mathcal{K}^{BC}_{+,1}(x;K,\rho)$ underlying the robust
bias-corrected estimator defined in Section S.5.3 of the supplement to \cite{calonico2020ej}. Replacing the
sample moment matrices in these weights with their limits
gives the equivalent kernel
corrected for bias corresponding to $\overline{\mathcal{K}}^{BC}_{+,1}(x;K,\rho)$. This representation provides a useful alternative interpretation of
our construction. The weights $W_{i,h}(x)$ and $W_{i,2,2,b}(x)$ are
the finite-sample local-linear and local-quadratic weighting components
associated with $\widehat m_{1,h}(x)$ and
$\widehat m_{2,b}^{(2)}(x)$, respectively, and  their robust combination in \eqref{rbc-weight} produces, up to a multiplicative factor, the finite-sample weights of the bias-corrected local-linear estimator used in the robust Wald procedure of \cite{calonico2020ej}.
Another related interpretation is provided by the influence-function
principle studied by \cite{bravo2020two}. The main insight is to correct the original estimating equations by using the influence function of the plug-in sample moment as the new moment function.
In the RDD setting, this
correction leads, up to normalization, to the same weighting structure
as the bias-corrected equivalent kernel
$\overline{\mathcal K}^{BC}_{+,1}(x;K,\rho)$.

The above equivalence has two important implications. First, although our robust weights are motivated by the resampling interpretation of EL, they also arise naturally from the equivalent-kernel representation of the underlying bias-corrected estimator and from the associated influence-function correction. This provides further theoretical support for the proposed weighting schemes and helps explain the practical effectiveness. Second, and more importantly, the resampling perspective gives a simple and informative construction of these weights and shows directly how bias-corrected EL ratios can be built within the empirical likelihood framework.
}
\end{remark}

\begin{remark}
{\color{black}
Although the above discussion focuses on local-linear estimation with
local-quadratic bias correction, the proposed framework naturally extends
 to higher-order local polynomial smoothers. Specifically, using local polynomial smoothers of
orders $p$ and $p+1$ for estimation and bias correction, respectively,
leads to robust weights
\begin{equation*}
\label{rbc-weight-p}
W_{i,p,h,b}^\star(x) = W_{i,0,p,h}(x) - W_{i,p+1,p+1,b}(x) \frac{1}{nb^{p+1}} \sum_{k=1}^n W_{k,0,p,h}(x) (X_k - x)^{p+1}.
\end{equation*}
These weights coincide with the finite-sample weights 
$\mathcal{K}^{BC}_{+,p}(x;K,\rho)$, as defined in Section S.5.3 of the supplement to \cite{calonico2020ej}.}
\end{remark}


\subsection{Robust empirical likelihood}
\label{sec:3.2}

Building on the robust weights constructed in Section~\ref{sec:3.1}, we propose two refined versions of the bias-corrected empirical log-likelihood ratios for  $m(x)$ in model (\ref{eq1}). 
To be specific, the {\bf T}aylor-expansion-based {\bf R}obust EL ratio and {\bf D}ifference-based {\bf R}obust EL ratio are respectively given by
\begin{equation*}
\begin{aligned}
l_{\TR}(\theta) = - 2 \max \bigg \{\sum_{i = 1}^n \log(np_i) \,\bigg|\, p_i \ge 0, \,\sum_{i=1}^n p_i = 1, \,\sum_{i=1}^n p_i W_{i,h,b}^\star(x) (Y_i - \theta) = 0 \bigg \},
\end{aligned}
\end{equation*}
\begin{equation*}
\begin{aligned}
l_{\DR}(\theta) = - 2 \max \bigg\{\sum_{i = 1}^n \log(np_i) \,\bigg|\, p_i \ge 0, \,\sum_{i=1}^n p_i = 1,\, \sum_{i=1}^n p_i  W_{i,h,b}^{\diamond}(x) (Y_i - \theta) = 0 \bigg\},
\end{aligned}
\end{equation*}
where  $ W_{i,h,b}^{\star}(x)$ is  defined in \eqref{rbc-weight} and $ W_{i,h,b}^{\diamond}(x)$ is specified in \eqref{rbc-weight2}.
Before presenting their limiting distributions, we impose the refined asymptotic framework formally in Condition~\ref{cond5}.
\begin{condition}\label{cond5}
The bandwidth $ h $ and the pilot bandwidth $ b $ satisfy
$$
\frac{h}{b} \to \kappa \in [0,1], \text{   as   }  n \rightarrow \infty.
$$
\end{condition}

\begin{theorem}
\label{th2}
Assume that Conditions {\rm\ref{cond1}--\ref{cond5}} hold, $m^{(2)}(x) \neq 0$ and $x$ is an interior point. Then, as $n \to \infty,$
\begin{equation*}
l_{{ \TR}}\left\{m(x)\right\} \stackrel{\mathcal{D}}{\rightarrow} \chi_{1}^{2}, \quad \text{and} \quad l_{\DR}\left\{m(x)\right\} \stackrel{\mathcal{D}}{\rightarrow} \chi_{1}^{2}.
\end{equation*}
\end{theorem}
Unlike the conventional bias-corrected  EL ratios in Section~\ref{sec:2.2}, Theorem \ref{th2} demonstrates that the proposed $l_{\TR}$ and $l_{\DR}$  asymptotically follow $\chi_1^2$
even 
when the bandwidth ratio $h/b$ converges to a positive constant, ensuring the validity of the proposed robust EL methods. 
This highlights that the robust weighting schemes  in \eqref{rbc-weight} and \eqref{rbc-weight2} can internally adjust for the additional variability from bias correction,   
thus allowing for simple-yet-robust EL formulations that preserve the pivotal asymptotic properties across a broader range of bandwidth choices. Supported by Theorem~\ref{th2}, we construct the robust EL confidence interval for $m(x)$ at the confidence level $1-\alpha$ as 
\[
I_{\alpha,\TR} = \big\{\theta : l_{\TR}(\theta) \leq \chi_1^2(\alpha) \big\},~~\mbox{or}~~
I_{\alpha,\DR} = \big\{\theta : l_{\DR}(\theta) \leq \chi_1^2(\alpha) \big\},
\]
where $\chi_1^2(\alpha)$ represents the upper $\alpha$-quantile of the $\chi_1^2$-distribution.

\begin{remark}
{\color{black}
For Wald-type inference on the RDD setting, Remark~7 of
\cite{Calonico2014} shows that, when $h=b$, a bias-corrected local
polynomial estimator of order $p$, together with its
bias-estimation-aware variance, is equivalent to the
corresponding local polynomial estimator of order $p+1$. This suggests
an EL implementation using higher-order local polynomial weights.
In particular, local-quadratic weights yield an EL ratio that retains
the Wilks phenomenon for bandwidths of order $n^{-1/5}$; see Remark~4
of \cite{ma2020empirical}. Our general framework includes this special case \(h=b\), where our robust weights (\ref{rbc-weight}) reduce to local-quadratic weights (up to a multiplicative factor), while allowing greater flexibility through \(h/b\to\kappa\in(0,1]\). Note that such  ratio has
important inferential implications. As shown by \cite{calonico2020ej}, the ratio
\(h/b\) affects the leading constants governing coverage accuracy and
interval length, and the length-optimal choice need not be one. Consistently, our simulations in Section~\ref{sec:sim} indicate that although \(h=b\) delivers coverage close to the nominal level, allowing \(b>h\) can achieve comparable coverage with shorter average confidence intervals. For practical applications, we recommend adopting our proposed robust framework. 
}
\end{remark}

\section{Regression discontinuity designs}
\label{sec:RDD}
We are now ready to apply our proposed EL-based inference framework to RDD analysis. Our goal is to develop a robust inference procedure for the regression discontinuity (RD) average treatment effect at the cutoff.
Following the methodological framework in Section~\ref{sec:methodology}, we first introduce new weights and robust EL ratios for two main RDD settings, i.e., sharp and fuzzy RDD, in Sections~\ref{sec:4.1} and~\ref{sec:4.2}, respectively. We then 
establish the limiting distributions and based on which construct the corresponding EL-based confidence intervals in Section~\ref{sec:4.3}.


\subsection{Sharp RDD}
\label{sec:4.1}

In the canonical sharp RDD, we consider a random sample $\big\{\big(Y_i(0), Y_i(1), X_i\big), i = 1, \ldots, n\big\}$ drawn from the joint distribution of the triplet $\big(Y(0), Y(1), X\big)$.  The covariate $X_i$, often referred to as the running or forcing variable, determines whether unit $i$ receives the treatment ($X_i \geq c$) or not ($X_i < c$). For simplicity, we set the cutoff point $c = 0$ without loss of generality.  The random variables $Y_i(1)$ and $Y_i(0)$ denote the potential outcomes for unit $i$ under treatment and control, respectively.
In practice, we observe the random sample $\big\{(Y_i, X_i), i = 1, \ldots, n\big\}$, where the observed outcome $Y_i$ for each unit $i$ is given by:
\[ Y_i = (1 - T_i) Y_i(0) + T_i Y_i(1), \]
with $T_i = I(X_i \geq 0)$ indicating treatment receipt and $I(\cdot)$ denoting the indicator function. Given that the conditional expectation functions $E\left\{Y_i(1)  \mid X = x\right\}$ and $E\left\{Y_i(0) \mid X = x\right\}$ are continuous at $x=0$, the sharp RD average treatment effect at the cutoff is identified as:
\[
\tau_{\SRD} \equiv E\left\{Y_i(1) - Y_i(0) \mid X = 0\right\} = \lim_{x \to 0^+} E(Y_i \mid X_i = x) - \lim_{x \to 0^-} E(Y_i \mid X_i = x).
\]
Write $\mu(x) = E(Y_i \mid X_i = x), \mu_+ = \lim_{x \to 0^+} \mu(x)$ and $ \mu_- = \lim_{x \to 0^-} \mu(x)$.
Then
$
\tau_{\SRD} = \mu_+ - \mu_-.
$
The problem can thus be reframed as inferring the difference between the regression functions of the outcome, given the running variable, at the cutoff for the control and treated groups.


We now introduce the local linear weights for treated and control units, respectively,
\begin{align*}
W_{i,h}^{+} &= I(X_i \geq 0) K_h(X_i) \bigg(S_{2,h}^{+} - S_{1,h}^{+}\frac{X_i}{h}\bigg), 
~~~S_{j,h}^{+} = \frac{1}{n} \sum_{i=1}^n I(X_i \geq 0) K_h(X_i)\left(\frac{X_i}{h}\right)^j, \\
W_{i,h}^{-} &= I(X_i < 0) K_h(X_i) \left(S_{2,h}^{-} - S_{1,h}^{-}\frac{X_i}{h}\right),
~~~S_{j,h}^{-} = \frac{1}{n} \sum_{i=1}^n I(X_i < 0) K_h(X_i)\left(\frac{X_i}{h}\right)^j.
\end{align*}
Let $\bZ_i(\theta,a) = \left(W_{i,h}^{+} \left(Y_i-\theta - a\right),~W_{i,h}^{-}\left(Y_i-a\right)\right)^\T.$ The original empirical log-likelihood function, introduced by \cite{Otsu2015}, is  given by
\begin{align}
\nonumber
l_\SRD(\theta,a) = - 2 \max \bigg\{\sum_{i=1}^n\log(np_i) \,\bigg|\, p_i \ge 0, \, \sum_{i=1}^n p_i = 1, \,\sum_{i=1}^n p_i \bZ_i(\theta,a) = 0 \bigg\}.
\end{align}
and the corresponding empirical log-likelihood ratio for $\tau_{\SRD}$ is defined as
\[
l_\SRD(\theta) = \inf_{a \in \mathcal{A}} l_\SRD(\theta,a),
\]
where $\mathcal{A}$ is the compact parameter space for $a$.

To implement the proposed robust EL approach, we employ a pilot bandwidth $b$
and introduce $W_{i,2,2,b}^{+}$ and $W_{i,2,2,b}^{-}$ as counterparts of $W_{i,2,2,b}(0)$, computed using the subsamples $\{i: X_i \ge 0\}$ and $\{i: X_i < 0\}$, respectively. Inspired by the robust weights introduced in Section~\ref{sec:3.1}, for $i=1,\ldots,n,$ we define
\[
W_{i,h,b}^{\star+} = W_{i,h}^{+} - n^{-1}b^{-2} W_{i,2,2,b}^{+}\sum_{k=1}^n W_{k,h}^{+}X_k^2 , ~~~ W_{i,h,b}^{\star-} = W_{i,h}^{-} - n^{-1}b^{-2}W_{i,2,2,b}^{-}\sum_{k=1}^n W_{k,h}^{-}X_k^2 ,
\]
and write
$$
\bZ_i^\star(\theta,a) = \left(W_{i,h,b}^{\star+} \left(Y_i-\theta - a\right),~W_{i,h,b}^{\star-}\left(Y_i-a\right)\right)^\T.
$$
We then propose the Taylor-expansion-based robust empirical log-likelihood function
\begin{equation}
\label{sharp-rdd-1}
\begin{aligned}
l_{\SRD,\TR}(\theta,a) = - 2 \max \bigg\{\sum_{i=1}^n\log(np_i) \,\bigg|\,  p_i \ge 0, \,\sum_{i=1}^n p_i = 1, \,\sum_{i=1}^n p_i \bZ_i^{\star}(\theta,a)=0\bigg\},
\end{aligned}
\end{equation} 
and the corresponding robust empirical log-likelihood ratio for $\tau_{\SRD}$ as
\begin{equation}
\label{sharp-rdd-2}
l_{\SRD,\TR}(\theta) = \inf_{a \in \mathcal{A}} l_{\SRD,\TR}(\theta,a).
\end{equation}
To compute (\ref{sharp-rdd-2}), we employ the Lagrange multiplier method, yielding
\begin{equation}
\label{sharp-rdd-com-1}
l_{\SRD,\TR}(\theta) = 2 \sum_{i=1}^n \log \big\{1 + \blambda^\T \bZ_{i}^{\star}(\theta,a)\big\},
\end{equation}
where $\blambda $ and $a$ are determined by 
\begin{equation}
\label{sharp-rdd-com-2}
\begin{aligned}
\sum_{i=1}^n \frac{\bZ_i^\star(\theta,a)} {1 +\blambda^\T\bZ_i^\star(\theta,a)} = 0 ~~~\mbox{and}~~~
\sum_{i=1}^n \frac{\blambda^\T \bW_i^\star} {1 +\blambda^\T\bZ_i^\star(\theta,a)} = 0, 
\end{aligned}
\end{equation}
with $\bW_i^\star = (W_{i,h,b}^{\star+}, W_{i,h,b}^{\star-})^\T.$
Consequently, these equations can be effectively solved using the Newton-Lagrange algorithm or nested algorithm, as discussed in \cite{Owen2001book}.

Similarly, let $W_{i,0,2,b}^{+}(x)$ and $W_{i,0,2,b}^{-}(x)$ be counterparts of $W_{i,0,2,b}(x)$ using the subsamples $\{i: T_i=1\}$ and $\{i: T_i = 0\}$, respectively, and write 
\begin{equation*}
\begin{aligned}
W_{i,h,b}^{\diamond+}  &= W_{i,h}^+ - n^{-1}\sum_{k=1}^n W_{k,h}^+ \big \{W_{i,0,2,b}^{+}(X_k) - W_{i,0,2,b}^{+}(0)\big\},\\
W_{i,h,b}^{\diamond-}  &= W_{i,h}^- - n^{-1}\sum_{k=1}^n W_{k,h}^- \big \{W_{i,0,2,b}^{-}(X_k) - W_{i,0,2,b}^{-}(0)\big\}.
\end{aligned}
\end{equation*}
The difference-based robust empirical likelihood ratio for $\tau_{\SRD}$, denoted as $l_{\SRD,\DR}(\theta)$, can be defined and computed analogously by replacing $W_{i,h,b}^{\star+}$ and $W_{i,h,b}^{\star-}$ in \eqref{sharp-rdd-1}--\eqref{sharp-rdd-com-2} with $W_{i,h,b}^{\diamond+}$ and $W_{i,h,b}^{\diamond-}$, respectively. 




\begin{remark}
   Given the representation of $W_{i,h,b}^{\star+}$, the bias-corrected local linear estimator $\widehat \tau_{\SRD}$ for $\tau_{\SRD}$ can be expressed as
$
\widehat \tau_{\SRD} = \widehat \mu_+ - \widehat \mu_-,
$
where the estimators $\widehat \mu_+$ and $\widehat \mu_-$ satisfy the following moment equations:
\begin{align*}
\sum_{i=1}^n W_{i,h,b}^{\star+}\left(Y_i-\widehat{\mu}_{+}\right) = 0, ~~
\sum_{i=1}^n W_{i,h,b}^{\star-}\left(Y_i-\widehat{\mu}_{-}\right) = 0.
\end{align*}
\cite{Calonico2014} derived the asymptotic variance of $\widehat \tau_{\SRD}$ and constructed the  normal-approximation-based robust confidence intervals for $\tau_{\SRD}$. \textcolor{black}{However, the variance introduced in their work is complex and requires additional estimation. \cite{calonico2020ej} later developed a fixed-$n$ variance estimator that simplifies computation and avoids the need for additional tuning parameters.
Our robust EL proposal instead provides a distinct route, in which both bias correction and the variability induced by bias estimation are incorporated directly into the EL formulation.
This leads to valid inference without a separate studentization step.}
\end{remark}

\subsection{Fuzzy RDD}
\label{sec:4.2}
Unlike the sharp RDD, which assumes that treatment status $T_i$  is fully assigned based on the forcing variable $X_i$, fuzzy RDD accommodates scenarios where $X_i$ influences but does not completely determine the treatment assignment. Specifically, in a fuzzy RDD, the conditional probability of receiving treatment changes discontinuously at the cutoff point $X_i = 0$:
\[
\lim_{x \to 0^+} P(T_i=1 \mid X_i=x) \neq \lim_{x \to 0^-} P(T_i=1 \mid X_i=x).
\]
To formally describe this setting,  we consider a random sample $\big(Y_i(1), Y_i(0), T_i(1), T_i(0), X_i\big)$ for $i = 1,\ldots,n$  from the joint distribution of $\big(Y(1), Y(0), T(1), T(0), X\big)$, where the treatment status of unit $i$
 is determined by
$
T_i = T_i(0)I(X_i < 0) + T_i(1)I(X_i \geq 0),
$
with $T_i(0), T_i(1) \in \{0,1\}$. The observed data consist of the sample $\{(Y_i, T_i, X_i), i = 1,\ldots,n\}$.
According to \cite{han2001}, under appropriate conditions, the average treatment effect at the cutoff in a fuzzy RDD is nonparametrically identified as
\[
\tau_{\FRD} = \frac{\lim_{x \to 0^+} E\left(Y_i \mid X_i=x\right) - \lim_{x \to 0^-} E\left(Y_i \mid X_i=x\right)}{\lim_{x \to 0^+} E\left(T_i \mid X_i=x\right) - \lim_{x \to 0^-} E\left(T_i \mid X_i=x\right)} 
\equiv \frac{\mu_{Y+} - \mu_{Y-}}{\mu_{T+} - \mu_{T-}}.
\]
The original empirical log-likelihood ratio, proposed by \cite{Otsu2015}, for $\tau_{\FRD}$ can then be  formulated as
\begin{equation*}
l_{\FRD}(\theta) = \inf_{(a,b_+,b_-) \in \mathcal{A} \times [0,1] \times [0,1]} l_{\FRD}(\theta,a,b_+,b_-),
\end{equation*}
where $l_{\FRD}(\theta,a,b_+,b_-)$ is defined as
\begin{equation}
\nonumber
\begin{aligned}
l_{\FRD}(\theta,a,b_+,b_-) = & - 2 \max \bigg\{\sum_{i=1}^n\log(np_i) \,\bigg|\, p_i \ge 0, \,\sum_{i=1}^n p_i = 1, \\
& \sum_{i=1}^n p_i W_{i,h}^{+} \left(Y_i-\theta b_+ - a\right) = 0, \,\sum_{i=1}^n p_i W_{i,h}^{+} \left(T_i- b_+\right) = 0,\\
& \sum_{i=1}^n p_i W_{i,h}^{-} \left(Y_i-\theta b_- - a\right)=0, \,\sum_{i=1}^n p_i W_{i,h}^{-} \left(T_i- b_-\right)=0 \bigg\}.
\end{aligned}
\end{equation}

Adopting the robust approach similar to that used in Section~\ref{sec:4.1}, we replace $W_{i,h}^{+}$ and $W_{i,h}^{-}$ with their robust counterparts $W_{i,h,b}^{\star+}$ and $W_{i,h,b}^{\star-}$ and introduce the Taylor-expansion-based robust empirical log-likelihood function:
\begin{equation}
\label{fuzzy-rdd-1}
\begin{aligned}
l_{\FRD,\TR}(\theta,a,b_+,b_-) = & - 2 \max \bigg\{\sum_{i=1}^n\log(np_i) \,\bigg|\, p_i \ge 0, \,\sum_{i=1}^n p_i = 1, \\
& \sum_{i=1}^n p_i W_{i,h,b}^{\star+} \left(Y_i-\theta b_+ - a\right) = 0, \,\sum_{i=1}^n p_i W_{i,h,b}^{\star+}\left(T_i- b_+\right) = 0,\\
& \sum_{i=1}^n p_i W_{i,h,b}^{\star-}\left(Y_i-\theta b_- - a\right)=0, \,\sum_{i=1}^n p_i W_{i,h,b}^{\star-}\left(T_i- b_-\right)=0 \bigg\}.
\end{aligned}
\end{equation}
The corresponding robust empirical log-likelihood ratio for $\tau_{\FRD}$ is then given by
\begin{equation}
\label{fuzzy-rdd-2}
l_{\FRD,\TR}(\theta) = \inf_{(a,b_+,b_-) \in \mathcal{A} \times [0,1] \times [0,1]} l_{\FRD,\TR}(\theta,a,b_+,b_-).
\end{equation}
Let
$
\widetilde \bZ_i^\star(\theta,a) = \left(W_{i,h,b}^{\star+} \left(Y_i-\theta T_i - a\right),~W_{i,h,b}^{\star-}\left(Y_i- \theta T_i - a\right)\right)^\T.
$
We can reformulate  \eqref{fuzzy-rdd-2} as
\begin{equation}
\label{fuzzy-rdd-com-1}
l_{\FRD,\TR}(\theta) = 2 \sum_{i=1}^n \log \big\{1 + \blambda^\T \widetilde \bZ_i^\star(\theta,a)\big\},
\end{equation}
where the parameters $\blambda$ and $a$  satisfy 
\begin{equation}
\label{fuzzy-rdd-com-2}
\begin{aligned}
\sum_{i=1}^n \frac{\widetilde \bZ_i^\star(\theta,a)} {1 +\blambda^\T \widetilde \bZ_i^\star(\theta,a)} = 0, ~~~\mbox{and}~~~
\sum_{i=1}^n \frac{\blambda^\T \bW_i^\star} {1 +\blambda^\T \widetilde \bZ_i^\star(\theta,a)} = 0, 
\end{aligned}
\end{equation}
As a result, the computational algorithms that have been developed for the context of sharp RDD can be seamlessly applied to this scenario. Finally, substituting $W_{i,h,b}^{\star+}$ and $W_{i,h,b}^{\star-}$ in \eqref{fuzzy-rdd-1}--\eqref{fuzzy-rdd-com-2} with $W_{i,h,b}^{\diamond+}$ and $W_{i,h,b}^{\diamond-}$, respectively, we can define and compute the difference-based robust empirical likelihood ratio $l_{\FRD,\DR}(\theta)$ for $\tau_{\FRD}$.

\begin{remark} 
To study $\tau_{\FRD}$,
\cite{Calonico2014} investigated the bias-corrected local linear estimator $\widehat \tau_{\FRD}$ for $\tau_{\FRD}$, defined as:
\[
\widehat \tau_{\FRD} = \frac{\widehat \mu_{Y+} - \widehat \mu_{Y-}}{\widehat \mu_{T+} - \widehat \mu_{T-}},
\]
where the estimators $\widehat \mu_{Y+},~\widehat \mu_{Y-}, \widehat \mu_{T+}$ and $\widehat \mu_{T-}$ satisfy the following equations:
\begin{equation*}
\begin{aligned}
& \sum_{i=1}^n W_{i,h,b}^{\star+}\left(Y_i-\widehat{\mu}_{Y+}\right)=0, \quad \sum_{i=1}^n W_{i,h,b}^{\star-}\left(Y_i-\widehat{\mu}_{Y-}\right)=0, \\
& \sum_{i=1}^n W_{i,h,b}^{\star+}\left(T_i-\widehat{\mu}_{T+}\right)=0, \quad \sum_{i=1}^n W_{i,h,b}^{\star-}\left(T_i-\widehat{\mu}_{T-}\right)=0.
\end{aligned}
\end{equation*} 
While the estimator $\widehat \tau_{\FRD}$ has an analytic form, its variance estimator relies on multiple Taylor expansions and approximations, making it more complex than the variance derivation and estimation required in the sharp RDD case.
This further demonstrates a key advantage of our robust EL method, as it offers a unified framework for both RDD settings. 
Following similar arguments, our method can potentially be adapted to other RDD settings, such as kink RDD, which targets discontinuities in the first derivative of the regression function at the cutoff. 
\end{remark}

\subsection{Asymptotic properties}
\label{sec:4.3}
We first impose some regularity conditions on the sharp RDD setting.
\begin{condition}
\label{cond_sharp-rdd-1}
For some $\rho_0>0,$ the following holds in the neighborhood $(-\rho_0,\rho_0)$ around the cutoff $x = 0:$
\begin{itemize}
\item[{\rm (i)}] $E(Y_i^4|X_i = x) < \infty;$ 
\item[{\rm (ii)}] The density function $f(x)$ of X is continuous and bounded away from zero; 
\item[{\rm  (iii)}] $\mu_{+}(x) = E\{Y_i(1)|X_i = x\}$ and  $\mu_{-}(x) = E\{Y_i(0)|X_i = x\}$ are thrice continuously differentiable;
\item[{\rm (iv)}] The conditional variance $\sigma^2(x) = \var(Y_i|X_i =x )$ is right and left continuous at the cutoff $x =0$ and bounded away from zero.
\end{itemize}
 \end{condition}

\begin{condition}
\label{cond_sharp-rdd-2}
$\mathcal{A}$ is compact and the true value $\mu_{-}$ is an interior point of $\mathcal{A}.$
 \end{condition}

We then give an additional condition to address the fuzzy RDD setting.
 \begin{condition}
 \label{cond_fuzzy-rdd-1}
From some $\rho_0>0,$ the following holds in the neighborhood $(-\rho_0,\rho_0)$ around the cutoff $x = 0.$
\begin{itemize}
\item[{\rm (i)}]  $\mu_{T+}(x) = E\{T_i(1)|X_i = x\}$ and  $\mu_{T-}(x) = E\{T_i(0)|X_i = x\}$ are thrice continuously differentiable;
\item[{\rm (ii)}] The conditional variance $\sigma_T^2(x) = \var(T_i|X_i =x )$ is right and left continuous at $x =0$ and bounded away from zero.
\end{itemize}
\end{condition}

The above conditions are standard in the literature, with analogous assumptions imposed in \cite{Calonico2014} and \cite{Otsu2015}. In particular, Conditions~\ref{cond_sharp-rdd-1}(iii) and \ref{cond_fuzzy-rdd-1}(i) establish conventional smoothness requirements on the regression functions. These conditions play a pivotal role in governing the leading bias terms of the estimators in both settings. Furthermore, Conditions~\ref{cond_sharp-rdd-1}(iv) and \ref{cond_fuzzy-rdd-1}(ii) place standard constraints on the conditional variance of the observed outcome and treatment, respectively, allowing for potential heterogeneity across the threshold.
\begin{theorem}
\label{th3}
Assume that Conditions {\rm \ref{cond3}--\ref{cond_sharp-rdd-2}} hold. Then, as $n \to \infty,$
\begin{equation*}
l_{\SRD,\TR}\left(\tau_{\SRD}\right) \stackrel{\mathcal{D}}{\rightarrow} \chi_{1}^{2},
~~\mbox{and}~~l_{\SRD,\DR}\left(\tau_{\SRD}\right) \stackrel{\mathcal{D}}{\rightarrow} \chi_{1}^{2}.
\end{equation*}
\end{theorem}

\begin{theorem}
\label{th4}
Assume that Conditions {\rm \ref{cond3}--\ref{cond_fuzzy-rdd-1}} hold. Then, as $n \to \infty,$
\begin{equation*}
l_{\FRD,\TR}\left(\tau_{\FRD}\right) \stackrel{\mathcal{D}}{\rightarrow} \chi_{1}^{2},
~~\mbox{and}~~ l_{\FRD,\DR}\left(\tau_{\FRD}\right) \stackrel{\mathcal{D}}{\rightarrow} \chi_{1}^{2}.
\end{equation*}
\end{theorem}

Theorems~\ref{th3} and \ref{th4} reveal that the proposed EL ratios for both sharp and fuzzy RDD settings converge asymptotically to $\chi_1^2$, which nicely justifies the construction of robust EL-based confidence intervals. 
For the sharp RDD setting, the confidence intervals for $\tau_{\SRD}$ at the confidence level $1-\alpha$ are respectively given by
\[
I_{\alpha,\TR}^{(\SRD)} = \left\{\theta : l_{\SRD,\TR}(\theta) \leq \chi_1^2(\alpha) \right\},~~\mbox{and}~~
I_{\alpha,\DR}^{(\SRD)} = \left\{\theta : l_{\SRD,\DR}(\theta) \leq \chi_1^2(\alpha)\right\}.
\]
For fuzzy RDD, the corresponding confidence intervals for $\tau_{\FRD}$ are defined as
\[
I_{\alpha,\TR}^{(\FRD)} = \left\{\theta : l_{\FRD,\TR}(\theta) \leq \chi_1^2(\alpha) \right\},~~\mbox{and}~~
I_{\alpha,\DR}^{(\FRD)} = \left\{\theta : l_{\FRD,\DR}(\theta) \leq \chi_1^2(\alpha) \right\}.
\]

\subsection{Second-order asymptotics and bandwidth selection} \label{sec.second}

{\color{black} 
We further examine the second-order properties of the proposed robust EL procedure. These properties play an important role in characterizing the coverage accuracy of the resulting robust EL confidence intervals and can guide bandwidth selection by minimizing coverage error.

To facilitate this, we provide an explicit higher-order representation of
$l_{\SRD,\TR}(\theta_0)$ in Section~S.4 of the supplementary material.
Under standard regularity conditions
and by employing arguments analogous to those in \cite{chen2002confidence} and \cite{Otsu2015}, the resulting Edgeworth expansion suggests that
\begin{align*}
P\Big\{\tau_{\SRD} \in I_{\alpha,\TR}^{(\SRD)}\Big\} - (1- \alpha) = O\{nh^7 + h^2 + (nh)^{-1}\}.
\end{align*}
While this represents a slight improvement over the original EL-based confidence interval developed in \cite{Otsu2015}, which exhibits a coverage error of $O\{nh^5 + h^2 + (nh)^{-1}\}$, the overall gain remains limited. In contrast, \cite{calonico2020ej} derived the coverage error expansion for the Wald-based robust confidence interval as:
\begin{align*}
P\Big\{ \tau_{\SRD} \in I_{\text{RBC}}(h)\Big\} - (1- \alpha)= O\{nh^7 + h^3 + (nh)^{-1}\},
\end{align*}
which achieves a faster convergence rate. Our careful investigation demonstrates that 
an appropriate modification of the proposed robust EL procedure can possibly match this superior coverage error rate of $O_P\{nh^7 + h^3 + (nh)^{-1}\}$; see further discussion in Section~\ref{sec:diss}. Establishing such a result, however, would fundamentally alter the methodological framework of the current study. 
We therefore leave this comprehensive analysis as well as the  coverage-error-optimal bandwidth selection for future research.

To aid the numerical studies, we here adopt the MSE-optimal bandwidth selector of \cite{Calonico2014} and the coverage-error-optimal bandwidth selector of \cite{calonico2020ej}, both implemented in the R package \verb|rdrobust|, as convenient benchmarks. We also conduct sensitivity analyses over a broad range of bandwidths to assess the robustness and practical usefulness of the proposed robust EL methods in the simulation studies.
}
\section{Simulation studies}
\label{sec:sim}

We conduct a series of simulations to illustrate the performance of our proposed robust EL inference methods.  Sections~\ref{sec:sim.NR} and \ref{sec:sim.RDD} examine inference scenarios for 
$m(x)$ and the sharp RD effect $\tau_{\SRD}$, respectively.
In each scenario, we generate random samples $(X_i, Y_i)$ for $i=1, \dots, n$ from  \eqref{eq1}, where the measurement locations \(X_i\) and errors $\varepsilon_i$  are sampled independently from \(2 \mathcal{B}(2,4) - 1\), with \(\mathcal{B}(\beta_1, \beta_2)\) denoting a beta distribution with shape parameters \(\beta_1\) and \(\beta_2\), and the normal distribution \(\mathcal{N}(0, 0.1295^2)\), respectively.
The exact forms of $m(x)$ will be specified in the respective sections. For simplicity, we use the Epanechnikov kernel to compute the empirical log-likelihood ratios.
 
\subsection{Nonparametric regression}
\label{sec:sim.NR}
\begin{figure}[htbp]
    \centering
    \begin{subfigure}{0.9\linewidth}
        \centering
        \begin{subfigure}{0.49\linewidth}
            \centering  \includegraphics[width=\linewidth]{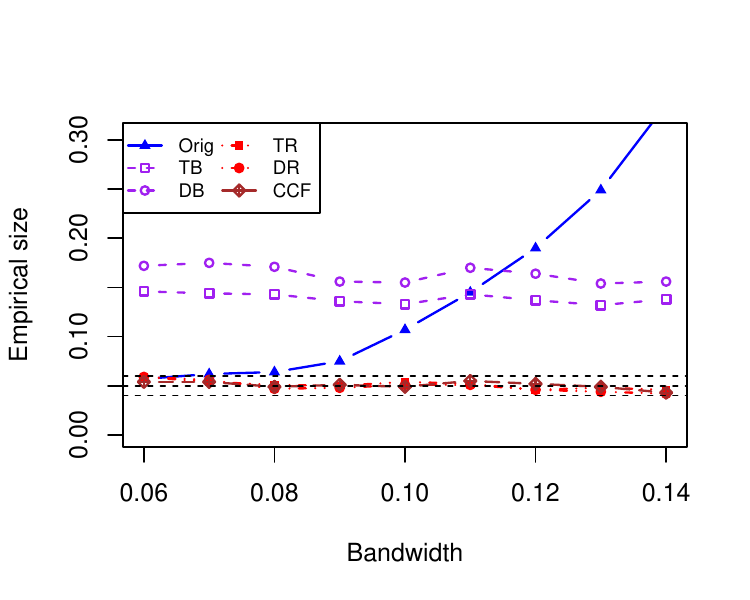}
        \end{subfigure}
        \hfill
        \begin{subfigure}{0.49\linewidth}
            \centering
            \includegraphics[width=\linewidth]{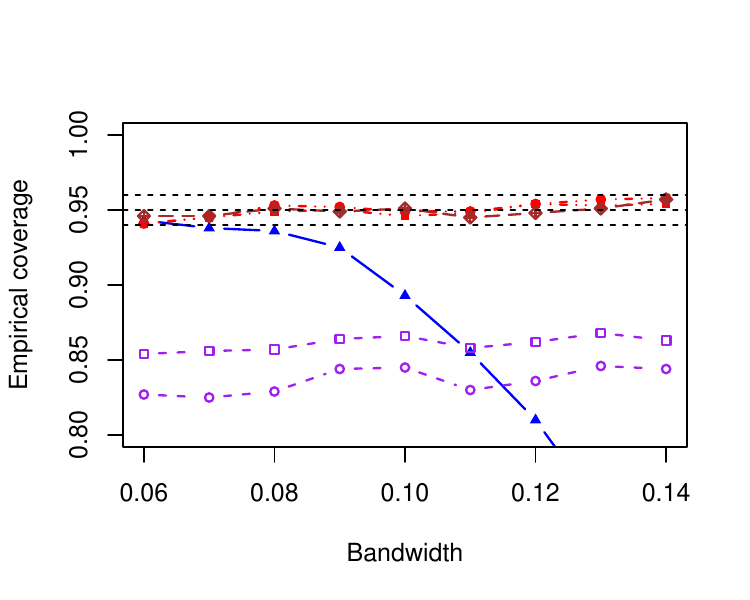}
        \end{subfigure}
        \vspace{-0.5cm}
        \caption{$n = 500, b = 1.2h$}
    \end{subfigure}

    \vspace{-0.8cm} 

    \begin{subfigure}{0.9\linewidth}
        \centering
        \begin{subfigure}{0.49\linewidth}
            \centering
            \includegraphics[width=\linewidth]{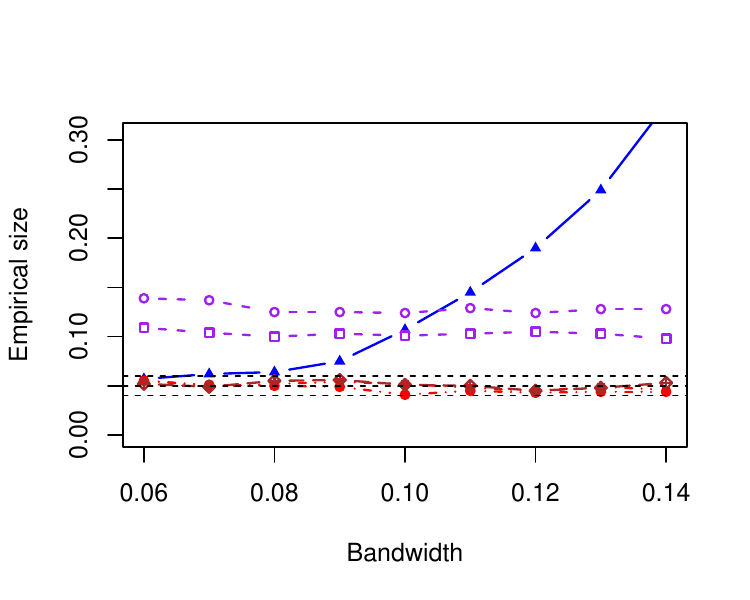}
        \end{subfigure}
        \hfill
        \begin{subfigure}{0.49\linewidth}
            \centering
            \includegraphics[width=\linewidth]{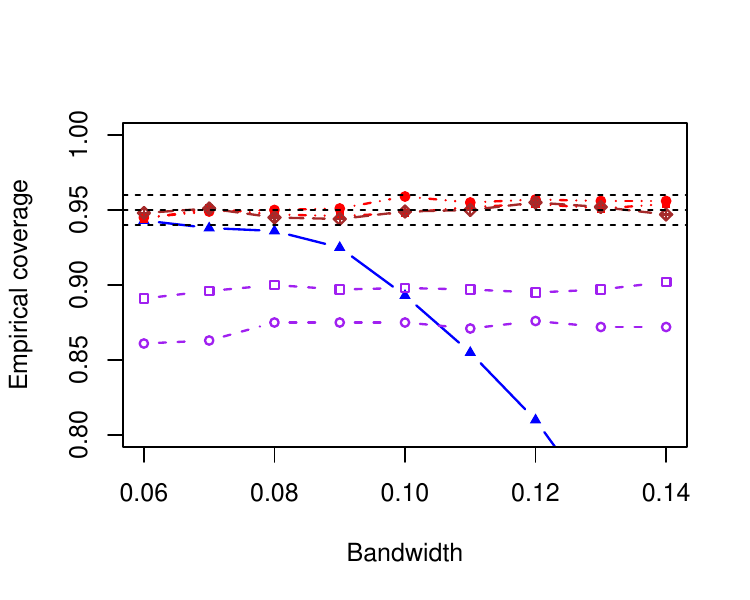}
        \end{subfigure}
        \vspace{-0.5cm}
        \caption{$n = 500, b = 1.5h$}
    \end{subfigure}
    
\vspace{-0.8cm}
 \begin{subfigure}{0.9\linewidth}
        \centering
        \begin{subfigure}{0.49\linewidth}
            \centering  \includegraphics[width=\linewidth]{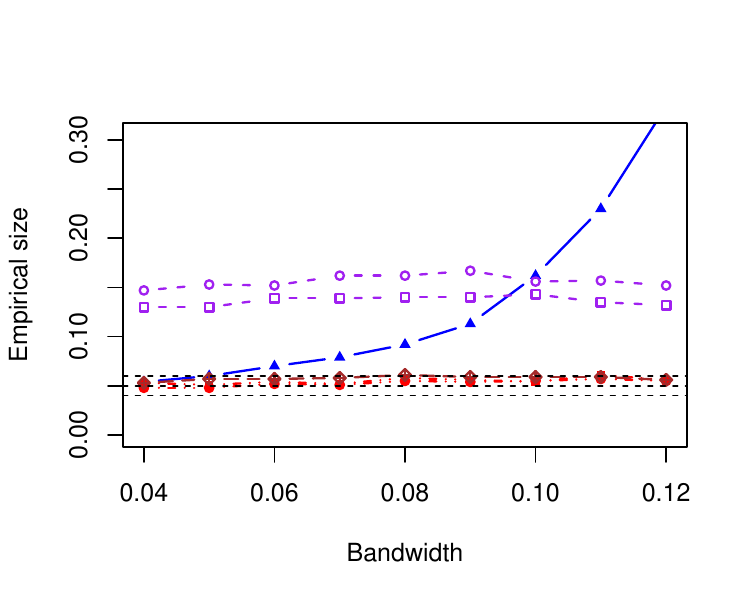}
        \end{subfigure}
        \hfill
        \begin{subfigure}{0.49\linewidth}
            \centering
            \includegraphics[width=\linewidth]{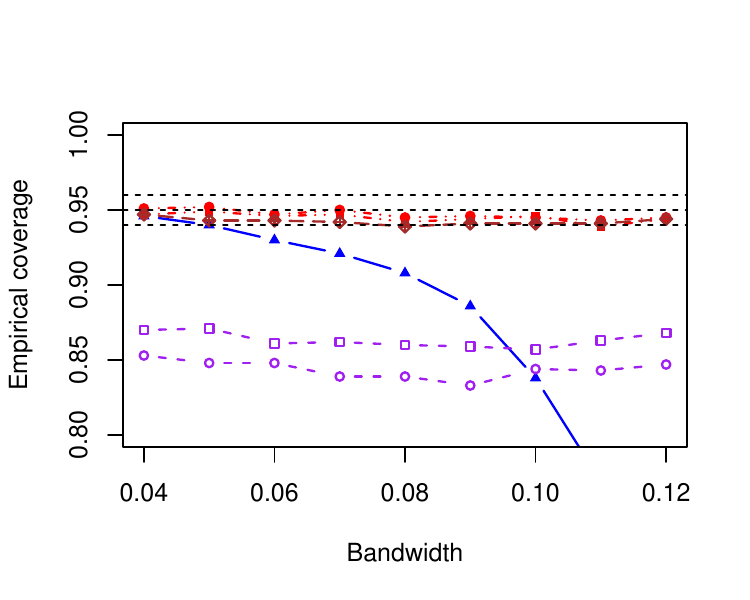}
        \end{subfigure}
        \vspace{-0.5cm}
        \caption{$n = 1000, b = 1.2h$}
    \end{subfigure}

\vspace{-0.8cm}
    \begin{subfigure}{0.9\linewidth}
        \centering
        \begin{subfigure}{0.49\linewidth}
            \centering
            \includegraphics[width=\linewidth]{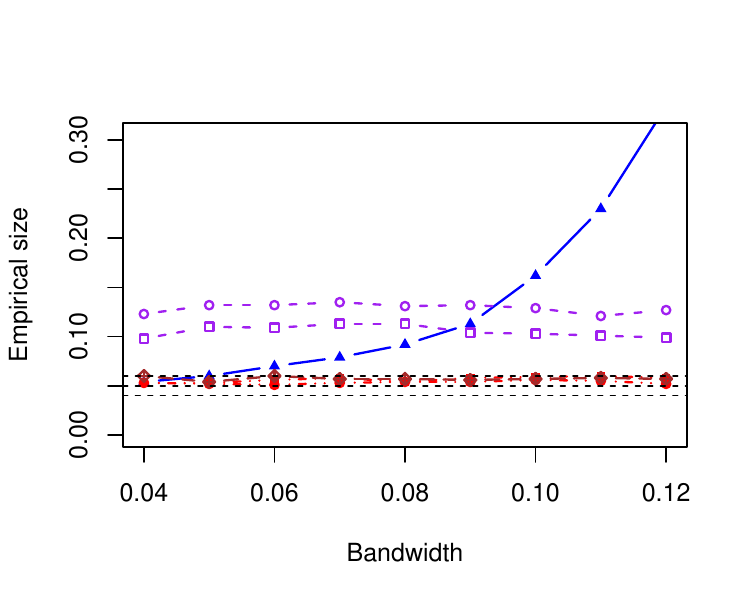}
        \end{subfigure}
        \hfill
        \begin{subfigure}{0.49\linewidth}
            \centering
            \includegraphics[width=\linewidth]{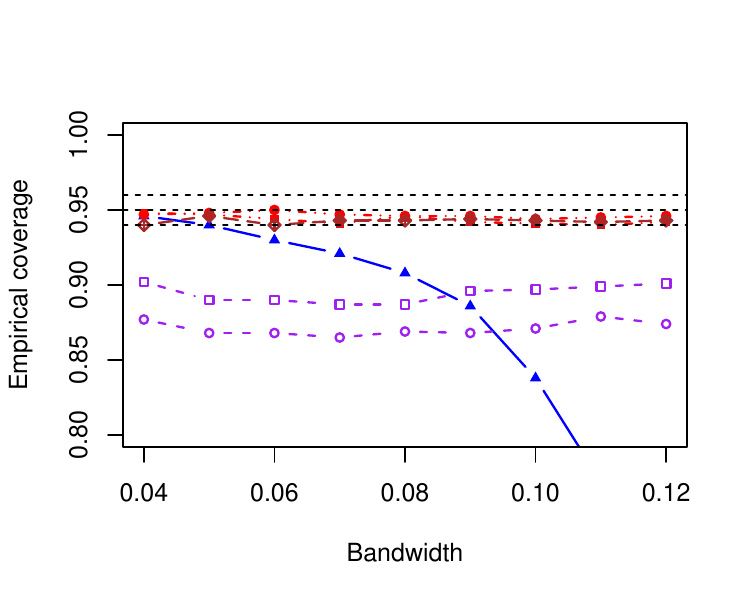}
        \end{subfigure}
        \vspace{-0.5cm}
        \caption{$n = 1000, b = 1.5h$}
    \end{subfigure}

    \caption{Plots of empirical sizes and empirical coverages as functions of bandwidth over 1000 simulation runs  at $x=-0.5$ for Model 1.}
    \label{fig:Case1}
\end{figure}

In this section, we consider the following two functional forms of $m(x)$, with a focus on inference at the point 
$x = -0.5$.
\newcounter{bean}
\setcounter{bean}{0}
\begin{center}
	\begin{list}
		{\textsc{Model} \arabic{bean}.}{\usecounter{bean}}
		\item \label{case1}  $m(x) = 0.25(x+1)^2-\sin(\pi x).$
		\item \label{case2} $m(x) = 0.3 \exp\big\{-4(2x+1)^2\big\} + 0.7 \exp\big\{-16(2x-1)^2\big\}.$ 

	\end{list}
\end{center} 
See also Figure~S1 of the supplementary material for visualizations of $m(x)$ and the evaluation points.   For each  model introduced above, we generate $n \in \{500,1000\}$ observations and replicate each simulation $1000$ times. We compare the performance of the proposed {\bf T}aylor-expansion-based {\bf R}obust EL method (TR) and {\bf D}ifference-based {\bf R}obust EL method (DR) against their conventional 
Bias-Corrected counterparts (denoted as TB and  DB, respectively),  the {\bf Orig}inal empirical likelihood method (Orig) introduced by \cite{chenQin2000},  \textcolor{black}{and the asymptotic robust bias correction approach
proposed by \cite{calonico2018jasa} (denoted as CCF).   The CCF method is implemented through R package \texttt{nprobust} \citep{nprobust}}.
To examine their robustness to the choice of bandwidths, we consider $h \in \{0.06, 0.08, \dots, 0.14\}$ for $n = 500$ and $h \in \{0.04, 0.06, \dots, 0.12\}$ for $n = 1000$, with three different pilot bandwidth settings  \textcolor{black}{$b =h$}, $1.2h$ and $1.5h$.
Figure~\ref{fig:Case1} reports empirical sizes for $m(x)$ at its true value under the $5\%$ nominal level and empirical interval coverages at the $95\%$ confidence level 
for five comparison methods with $b =1.2h$ and $1.5h$. Table~\ref{CI.case1} provides the average empirical interval lengths along with the empirical interval coverages. Since the results for Models~\ref{case1} and \ref{case2} exhibit similar trends, we only present numerical results for Model~\ref{case1} here and defer those for Model~\ref{case2} to Figure~S2 and Table~S1 in the supplementary material.
 \textcolor{black}{To complete the analysis, we also report Table~S2 in the supplementary material for a design with an asymmetric error distribution, where the errors \(\varepsilon_i\) are sampled independently from a scaled and standardized Gamma distribution. See Section~S.2 of the supplementary material for details.}

Several conclusions can be drawn from Figures~\ref{fig:Case1} and S2 and Tables~\ref{CI.case1} and S1. 
First, our proposed two robust methods (TR and DR), {\color{black}together with the asymptotic robust method (CCF),} perform equally well, achieving the most accurate empirical sizes and coverages across all settings. 
This demonstrates the effectiveness of our approaches and their robustness to both $h$ and $b$. 
Second, all non-robust competing methods suffer from significant size and coverage distortions. As expected, the Orig method is highly sensitive to bandwidth selection, with its performance deteriorating severely when $h$ becomes slightly larger. 
Meanwhile, both conventional methods (TB and DB) fail completely, as reflected by elevated empirical sizes and reduced empirical coverages across all values of $h$.
This pattern provides strong validation for Theorem~\ref{th1}, which states that $l_{\text{TB}}\{m(x)\}$ and $l_{\text{DB}}\{m(x)\}$ deviate from  $\chi_1^2$ when  the condition $h/b \to 0$ is violated. Interestingly, their performance improves slightly when $b = 1.5h$
 compared to $b = h$ and $1.2h$, likely due to the smaller ratio $h/b$ in the former case.
Thirdly, the average interval lengths decrease as $h$ and $n$ increase for all methods. This highlights the advantage of our proposed robust methods, which can produce narrower confidence intervals while maintaining correct size and coverage for a given $n$, as evidenced by the bolded results. In contrast, the Orig method, constrained by the requirement for undersmoothing,  results in wider intervals. 
 \textcolor{black}{Lastly, as also reflected in the bolded results, the average interval lengths of all robust methods (TR, DR, and CCT) further decrease as the ratio \(b/h\) increases, with their coverages remaining close to the nominal level. This suggests that allowing \(b>h\) has meaningful finite-sample implications and again supports the refined asymptotic framework in \eqref{eq.ratio}.}

\begin{table}[htbp]
	\caption{\label{CI.case1}  \textcolor{black}{Average interval lengths and empirical coverages (in parentheses) of $95\%$ confidence intervals over 1000 simulation runs  at $x=-0.5$ for Model 1. The best performances, with empirical coverages close to 95\% and short interval lengths, are in bold font.}
    }
	\begin{center}
	\vspace{-0.5cm}
		\resizebox{6in}{!}{

\begin{tabular}{ccccccccccccc}
\hline
\textbf{}                                                & \textbf{} & \multicolumn{5}{c}{{$n = 500$}}                                       & \textbf{} & \multicolumn{5}{c}{{$n = 1000$}}                                      \\ \cline{3-7} \cline{9-13} 
$h$                                                      & \textbf{} & {0.06} & {0.08} & {0.1} & {0.12} & {0.14} & \textbf{} & {0.04} & {0.06} & {0.08} & {0.1} & {0.12} \\
\hline
Orig                                                     &           & \textbf{0.071}        & 0.062         & 0.055        & 0.051         & 0.047         &           & \textbf{0.062}        & 0.050         & 0.043         & 0.039        & 0.036         \\
                                                         &           & \textbf{(0.943)}      & (0.936)       & (0.893)      & (0.810)       & (0.663)       &           & \textbf{(0.946)}       & (0.930)       & (0.908)       & (0.838)      & (0.666)       \\
\specialrule{0pt}{0.3em}{0.em}     &           & \multicolumn{11}{c}{ \textcolor{black}{Setting: $b = h$}}                                                                                                                                   \\
\specialrule{0pt}{0.3em}{0.em}   TB &           & 0.071         & 0.061         & 0.055        & 0.050         & 0.047         &           & 0.061         & 0.050         & 0.043         & 0.039        & 0.035         \\
                                                         &           & (0.815)       & (0.812)       & (0.811)      & (0.823)       & (0.813)       &           & (0.827)       & (0.825)       & (0.820)       & (0.815)      & (0.828)       \\
DB                                                       &           & 0.070         & 0.061         & 0.055        & 0.050         & 0.046         &           & 0.061         & 0.050         & 0.043         & 0.039        & 0.035         \\
                                                         &           & (0.801)       & (0.800)       & (0.807)      & (0.816)       & (0.811)       &           & (0.819)       & (0.814)       & (0.811)       & (0.806)      & (0.821)       \\
TR                                                       &           & 0.103         & 0.089         & 0.080        & 0.073         & \textbf{0.068}        &           & 0.089         & 0.073         & 0.063         & 0.056        & \textbf{0.051}        \\
                                                         &           & (0.936)       & (0.949)       & (0.948)      & (0.947)       & \textbf{(0.951)}       &           & (0.946)       & (0.951)       & (0.947)       & (0.944)      & \textbf{(0.947)}       \\
DR                                                       &           & 0.105         & 0.091         & 0.081        & 0.074         & \textbf{0.069}        &           & 0.090         & 0.074         & 0.064         & 0.057        & \textbf{0.052}        \\
                                                         &           & (0.938)       & (0.947)       & (0.949)      & (0.948)       & \textbf{(0.954)}       &           & (0.947)       & (0.949)       & (0.947)       & (0.941)      & \textbf{(0.948)}       \\
CCF                                                     &           & 0.102         & 0.089         & 0.079        & 0.072         & \textbf{0.067}        &           & 0.088         & 0.072         & 0.062         & 0.056        & \textbf{0.051}        \\
                                                         &           & (0.941)       & (0.945)       & (0.951)      & (0.946)       & \textbf{(0.946)}       &           & (0.938)       & (0.944)       & (0.946)       & (0.943)      & \textbf{(0.942)}       \\ 
\specialrule{0pt}{0.3em}{0.em}     &           & \multicolumn{11}{c}{Setting: $b = 1.2h$}                                                                                                                                \\
\specialrule{0pt}{0.3em}{0.em}   TB &           & 0.071         & 0.061         & 0.055        & 0.050         & 0.047         &           & 0.061         & 0.050         & 0.043         & 0.039        & 0.035         \\
                                                         &           & (0.854)       & (0.857)       & (0.866)      & (0.862)       & (0.863)       &           & (0.870)       & (0.861)       & (0.860)       & (0.857)      & (0.868)       \\
DB                                                       &           & 0.071         & 0.061         & 0.055        & 0.050         & 0.046         &           & 0.061         & 0.050         & 0.043         & 0.039        & 0.035         \\
                                                         &           & (0.827)       & (0.829)       & (0.845)      & (0.836)       & (0.844)       &           & (0.853)       & (0.848)       & (0.839)       & (0.844)      & (0.847)       \\
TR                                                       &           & 0.094         & 0.081         & 0.072        & 0.066         & \textbf{0.062}         &           & 0.081         & 0.066         & 0.057         & 0.051        & \textbf{0.047}         \\
                                                         &           & (0.941)       & (0.949)       & (0.946)      & (0.954)       & \textbf{(0.954)}       &           & (0.947)       & (0.946)       & (0.942)       & (0.946)      & \textbf{(0.944)}       \\
DR                                                       &           & 0.099         & 0.085         & 0.076        & 0.070         & \textbf{0.065}         &           & 0.085         & 0.069         & 0.060         & 0.054        & \textbf{0.049}         \\
                                                         &           & (0.941)       & (0.953)       & (0.949)      & (0.954)       & \textbf{(0.958)}       &           & (0.951)       & (0.947)       & (0.945)       & (0.945)      & \textbf{(0.945)}       \\
CCF                                                     &           & 0.093         & 0.080         & 0.072        & 0.066         & \textbf{0.061}        &           & 0.080         & 0.065         & 0.056         & 0.050        & \textbf{0.046}         \\
                                                         &           & (0.946)       & (0.951)       & (0.951)      & (0.948)       & \textbf{(0.957)}       &           & (0.947)       & (0.943)       & (0.939)       & (0.941)      & \textbf{(0.944)}       \\ 
\specialrule{0pt}{0.3em}{0.em}     &           & \multicolumn{11}{c}{Setting: $b =1.5 h$}                                                                                                                                \\
\specialrule{0pt}{0.3em}{0.em}   TB &           & 0.071         & 0.062         & 0.055        & 0.050         & 0.047         &           & 0.062         & 0.050         & 0.043         & 0.039        & 0.036         \\
                                                         &           & (0.891)       & (0.900)       & (0.898)      & (0.895)       & (0.902)       &           & (0.902)       & (0.890)       & (0.887)       & (0.897)      & (0.901)       \\
DB                                                       &           & 0.071         & 0.061         & 0.055        & 0.050         & 0.046         &           & 0.061         & 0.050         & 0.043         & 0.039        & 0.035         \\
                                                         &           & (0.861)       & (0.875)       & (0.875)      & (0.876)       & (0.872)       &           & (0.877)       & (0.868)       & (0.869)       & (0.871)      & (0.874)       \\
TR                                                       &           & 0.086         & 0.074         & 0.066        & 0.061         & \textbf{0.057}         &           & 0.074         & 0.060         & 0.052         & 0.047        & \textbf{0.043}         \\
                                                         &           & (0.944)       & (0.947)       & (0.948)      & (0.954)       & \textbf{(0.954)}      &           & (0.948)       & (0.944)       & (0.946)       & (0.941)      & \textbf{(0.942)}       \\
DR                                                       &           & 0.091         & 0.079         & 0.071        & 0.065         & \textbf{0.061}         &           & 0.079         & 0.064         & 0.055         & 0.050        & \textbf{0.046}        \\
                                                         &           & (0.945)       & (0.950)       & (0.959)      & (0.957)       & \textbf{(0.956)}       &           & (0.947)       & (0.950)       & (0.946)       & (0.944)      & \textbf{(0.946)}       \\
CCF                                                     &           & 0.085         & 0.074         & 0.066        & 0.060         & \textbf{0.056}        &           & 0.073         & 0.060         & 0.052         & 0.046        & \textbf{0.042}        \\
                                                         &           & (0.948)       & (0.945)       & (0.949)      & (0.955)       & \textbf{(0.947)}      &           & (0.940)       & (0.940)       & (0.943)       & (0.943)      & \textbf{(0.943)}   \\ \hline   
\end{tabular}
		}	
	\end{center}
\end{table}

\subsection{Regression discontinuity design}
\label{sec:sim.RDD}
We now turn to the sharp RDD setting and investigate Model~\ref{case3} below,  constructed using data in \cite{Lee2008} with  \textcolor{black}{ $\tau_{\SRD} = 0.04$}, at the boundary point $x = 0$. A similar model can also be found in \cite{Imbens2012} and \cite{Calonico2014}.
Further analysis of a simple low-order polynomial function, termed Model~4, is provided in Section~S.2 of the supplementary material. See Figure~S3 for the visualizations of $m(x)$ for both models.
\begin{center}
    \begin{list}{\textsc{Model} \arabic{bean}.}{\usecounter{bean}}
        \addtocounter{bean}{2} 
        \item \label{case3} \begin{singlespace}
 \begin{align*}
 m(x) =
\begin{cases} 
  \textcolor{black}{0.48} + 1.27x + 7.18x^2 + 20.21x^3 + 21.54x^4 + 7.33x^5 & \text{if } x < 0, \\
  \textcolor{black}{0.52} + 0.84x - 3.00x^2 + 7.99x^3 - 9.01x^4 + 3.56x^5 & \text{if } x \geq 0.
\end{cases}
 \end{align*}  
 \end{singlespace}
    \end{list}
\end{center}

\begin{figure}[htbp]
    \centering
    \begin{subfigure}{0.9\linewidth}
        \centering
        \begin{subfigure}{0.49\linewidth}
            \centering  \includegraphics[width=\linewidth]{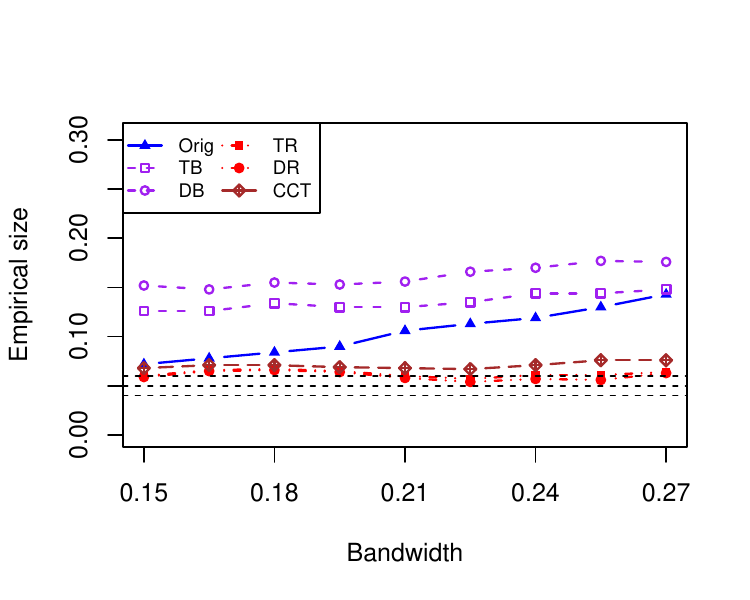}
        \end{subfigure}
        \hfill
        \begin{subfigure}{0.49\linewidth}
            \centering
            \includegraphics[width=\linewidth]{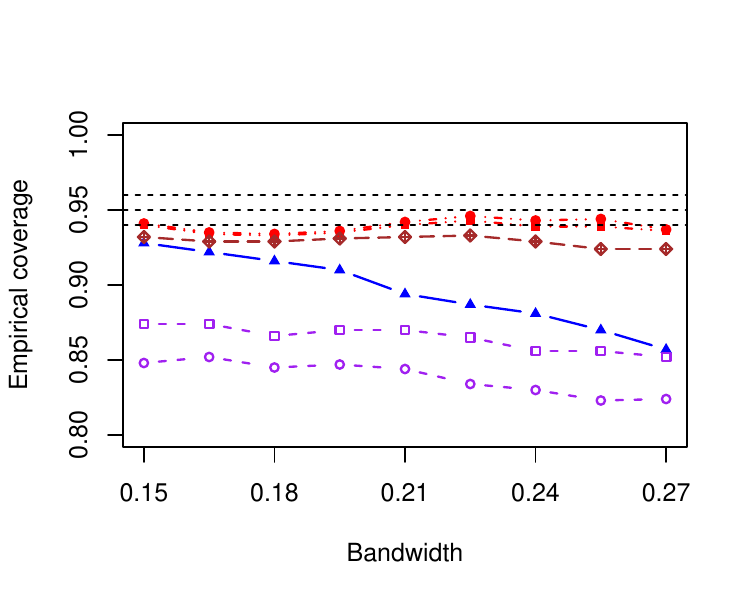}
        \end{subfigure}
        \vspace{-0.5cm}
        \caption{$n = 500, b = 1.2h$}
    \end{subfigure}

    \vspace{-0.8cm} 

    \begin{subfigure}{0.9\linewidth}
        \centering
        \begin{subfigure}{0.49\linewidth}
            \centering
            \includegraphics[width=\linewidth]{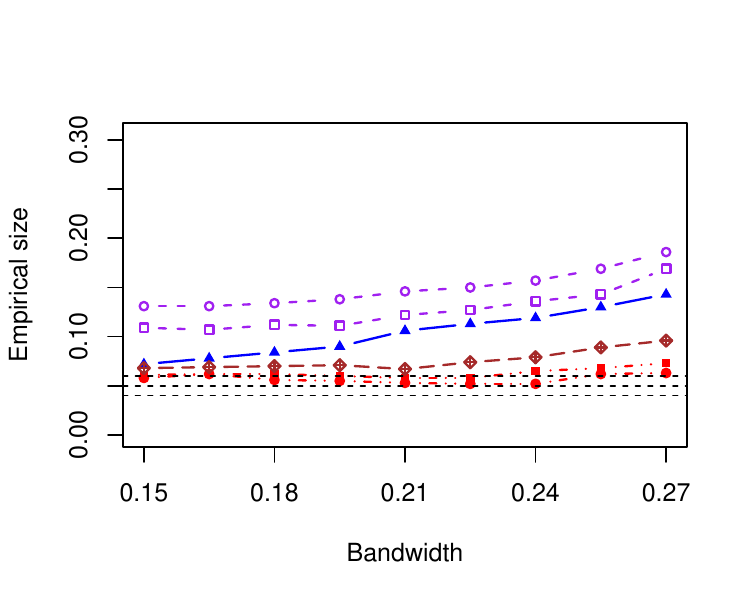}
        \end{subfigure}
        \hfill
        \begin{subfigure}{0.49\linewidth}
            \centering
            \includegraphics[width=\linewidth]{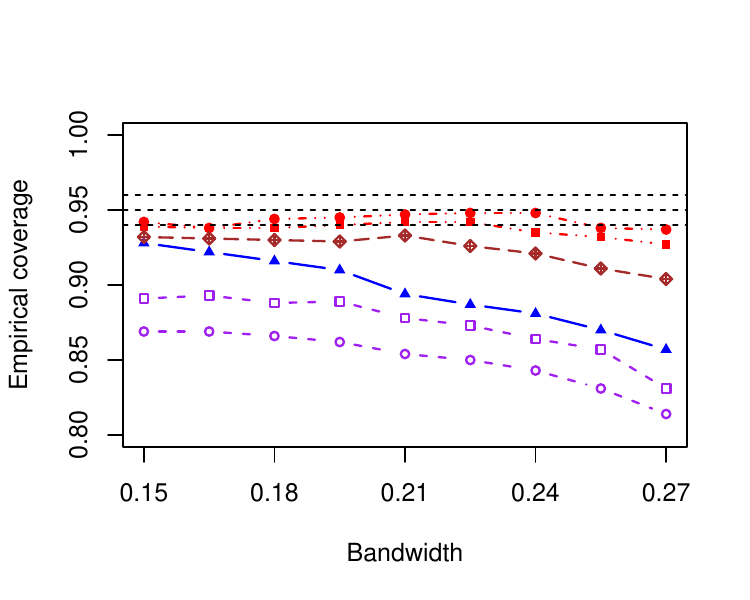}
        \end{subfigure}
        \vspace{-0.5cm}
        \caption{$n = 500, b = 1.5h$}
    \end{subfigure}
    
\vspace{-0.8cm}
 \begin{subfigure}{0.9\linewidth}
        \centering
        \begin{subfigure}{0.49\linewidth}
            \centering  \includegraphics[width=\linewidth]{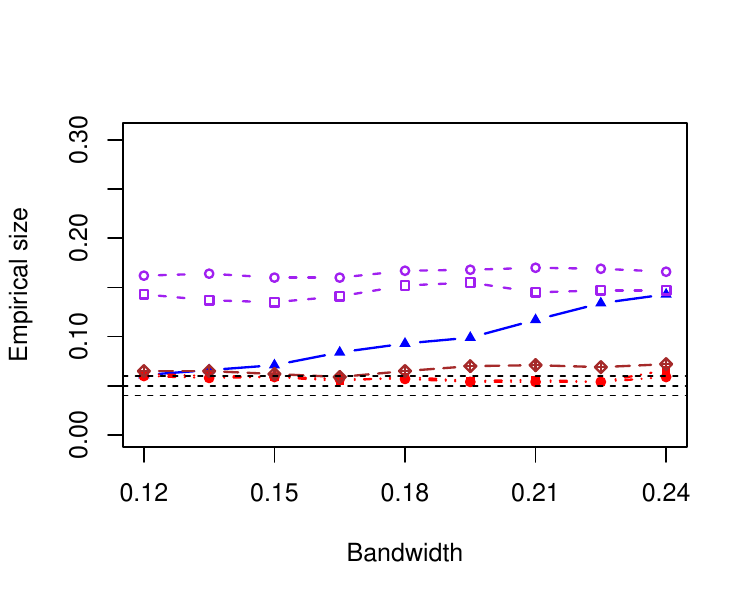}
        \end{subfigure}
        \hfill
        \begin{subfigure}{0.49\linewidth}
            \centering
            \includegraphics[width=\linewidth]{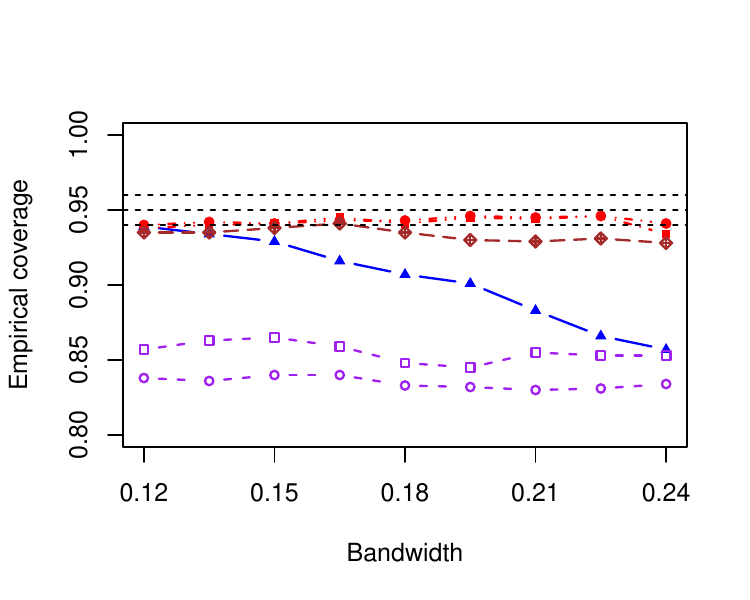}
        \end{subfigure}
        \vspace{-0.5cm}
        \caption{$n = 1000, b = 1.2h$}
    \end{subfigure}

\vspace{-0.8cm}
    \begin{subfigure}{0.9\linewidth}
        \centering
        \begin{subfigure}{0.49\linewidth}
            \centering
            \includegraphics[width=\linewidth]{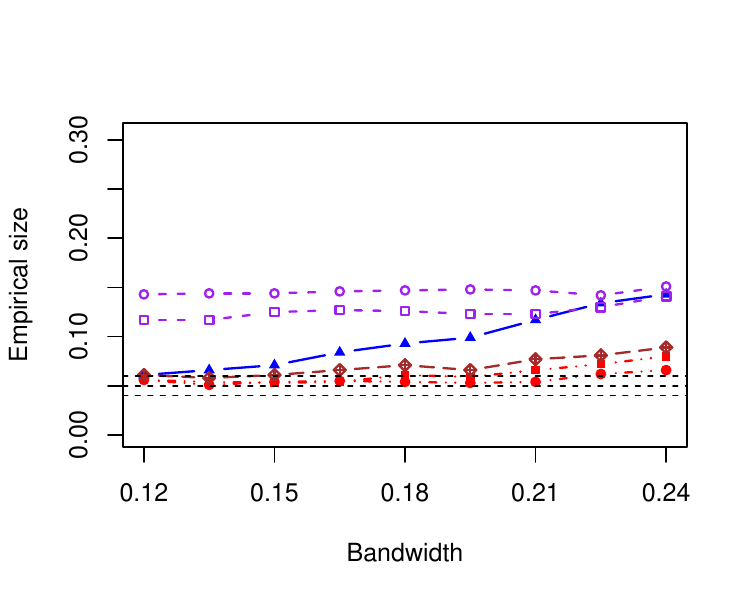}
        \end{subfigure}
        \hfill
        \begin{subfigure}{0.49\linewidth}
            \centering
            \includegraphics[width=\linewidth]{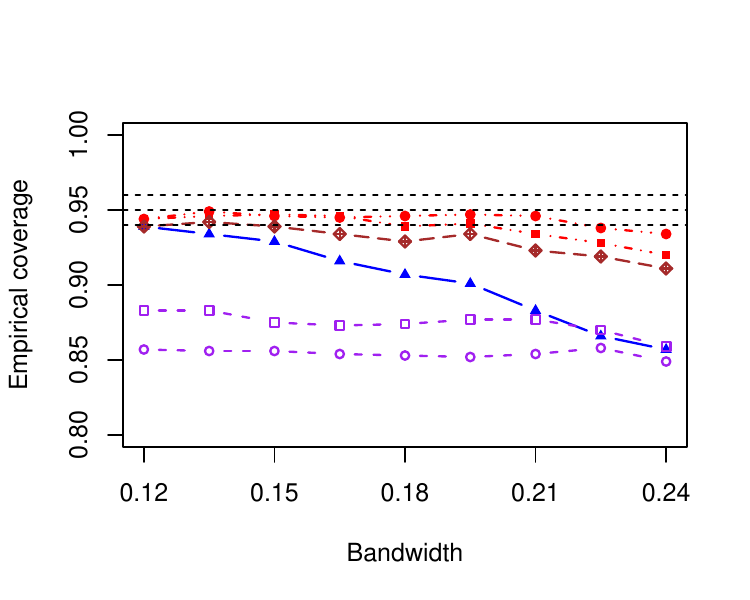}
        \end{subfigure}
        \vspace{-0.5cm}
        \caption{$n = 1000, b = 1.5h$}
    \end{subfigure}

    \caption{Plots of empirical sizes and coverages as functions of bandwidth over 1000 simulation runs for Model 3.}
    \label{fig:Case3}
\end{figure}

Similar to Section~\ref{sec:sim.NR}, we evaluate the performance of Orig, TB, DB, TR and DR methods under this sharp RDD setting.  \textcolor{black}{ For comparison, we also include the robust bias-corrected method of \cite{Calonico2014} (denoted by CCT). The EL-based methods are computed via Algorithm~1 in the supplementary material, and CCT is implemented using the \texttt{rdrobust} R package \citep{rdrobust}}.
To  assess the sensitivity of the competing methods to the bandwidth choices, 
Figure~\ref{fig:Case3} displays the  empirical sizes 
and  
empirical coverages 
based on 1000 simulations
over $h \in \{0.15,0.18, \dots, 0.27\}$ for $n =500$ and $h \in \{0.12,0.18, \dots, 0.24\}$ for $n = 1000$ under the settings $b = 1.2h$ and $1.5h$. 
Table~\ref{CI.case3} presents the average interval lengths and coverages  \textcolor{black}{with $b = h$}, $1.2h$ and $1.5h$.  The results align with the findings in Section~\ref{sec:sim.NR}. 
In particular, we observe a similar size distortion effect for the Orig method  as in \cite{Otsu2015}.
Overall, TR and DR outperform CCT across the range of $h$ values, 
delivering empirical sizes and coverages closer to the target levels while keeping comparable or shorter average interval lengths. For instance, when \( (n, h, b) = (1000, 0.21, 1.2h) \), our proposed methods (TR and DR) achieve more than \( 94.4\% \) empirical coverage with interval lengths of 0.188 and 0.190, respectively,  while CCT provides \( 93.5\% \) coverage with an interval length of 0.194 under \( (n, h, b) = (1000, 0.18, 1.2h) \). This pattern holds across other \( n \) and \( b \) settings. 

 \textcolor{black}{Table~\ref{CI.case3.time} reports the average computing time, in seconds, over 1000 simulation runs with \(b=1.2h\), on an AMD EPYC 9R14 processor at 3.69GHz, for both the test of \(\tau_{\SRD}=0.04\) and the construction of the 95\% confidence interval.
The EL-based procedures are naturally more computationally demanding than CCT, as interval construction requires numerical inversion of the EL ratio. Nevertheless, the computing times remain small with all EL-based confidence intervals constructed in less than 0.5 seconds. This indicates that Algorithm~1 provides a practically fast implementation of the EL-based methods for moderate sample sizes. As expected, the difference-based methods (DB and DR) require more computing time than the other EL-based methods, due to the additional computational cost in the construction of their difference-based weights.}

We further adopt the MSE-optimal bandwidth selectors for $h$ and $b$ in \cite{Calonico2014} to investigate the effectiveness of the various inference methods under such data-driven bandwidths, denoted as $\hat h_{\text{mse}}$ and $\hat b_{\text{mse}}$, respectively.   Table~\ref{opt.case3} provides numerical summaries for four different settings with $(h,b) = (\hat h_{\text{mse}}, \hat b_{\text{mse}}),$ \textcolor{black}{$ (\hat h_{\text{mse}}, \hat h_{\text{mse}}), $} $(\hat h_{\text{mse}}, 1.2 \hat h_{\text{mse}})$ and $ (\hat h_{\text{mse}}, 1.5\hat h_{\text{mse}})$ based on  10000 replications, including the empirical sizes for $\tau_{\SRD} = 0.04$ at $5\%$ nominal level, empirical coverages at $95\%$ confidence level and the average interval lengths.
Also reported are the average bandwidth values for the selected $h$ and $b$ in each setting. A few trends are apparent.  First, the Orig, TB and DB methods tend to significantly over-reject the null hypothesis, resulting in notable size and coverage errors. Second, the three robust methods demonstrate clear improvements in empirical sizes and coverages. Among them, DR performs the best, closely followed by TR, while CCT exhibits relatively inferior performance. This reaffirms the robustness and effectiveness of our proposed method in enhancing inference performance.  
 \textcolor{black}{We also present the results based on the coverage-error-optimal bandwidth selector of \cite{calonico2020ej} in Table~S3 of the supplementary material, where similar patterns are observed. Both the MSE-optimal and coverage-error-optimal bandwidth selectors are computed using \texttt{rdrobust}.}

\begin{table}[htbp]
	\caption{\label{CI.case3} \textcolor{black}{Average interval lengths and empirical coverages (in parentheses) of $95\%$ confidence intervals over 1000 simulation runs for Model 3.  The best performances, with empirical coverages close to 95\% and short interval lengths, are in bold font.}}
	\begin{center}
	\vspace{-0.5cm}
		\resizebox{6in}{!}{
       
\begin{tabular}{ccccccccccccc}
\hline
\textbf{} & \textbf{} & \multicolumn{5}{c}{$n = 500$}                   &  & \multicolumn{5}{c}{$n = 1000$}                  \\ \cline{3-7} \cline{9-13} 
$h$       & \textbf{} & 0.15    & 0.18    & 0.21    & 0.24    & 0.27    &  & 0.12    & 0.15    & 0.18    & 0.21    & 0.24    \\ \hline
Orig      &           & \textbf{0.227}   & 0.208   & 0.194   & 0.182   & 0.173   &  & \textbf{0.180}   & 0.161   & 0.147   & 0.137   & 0.129   \\
          &           & \textbf{(0.928)} & (0.916) & (0.894) & (0.881) & (0.857) &  & \textbf{(0.939)} & (0.929) & (0.907) & (0.883) & (0.857) \\
   \specialrule{0pt}{0.3em}{0.em}       &           & \multicolumn{11}{c}{\textcolor{black}{Setting: $b = h$}}                                                                \\
\specialrule{0pt}{0.3em}{0.em} TB        &           & 0.265   & 0.237   & 0.218   & 0.203   & 0.191   &  & 0.198   & 0.175   & 0.159   & 0.148   & 0.138   \\
          &           & (0.856) & (0.849) & (0.846) & (0.848) & (0.840) &  & (0.835) & (0.836) & (0.834) & (0.834) & (0.826) \\
DB        &           & 0.219   & 0.201   & 0.187   & 0.176   & 0.167   &  & 0.175   & 0.157   & 0.144   & 0.134   & 0.125   \\
          &           & (0.829) & (0.829) & (0.828) & (0.821) & (0.803) &  & (0.824) & (0.827) & (0.823) & (0.820) & (0.814) \\
TR        &           & 0.380   & 0.328   & 0.298   & 0.276   & \textbf{0.260}   &  & 0.282   & 0.246   & 0.223   & 0.206   & \textbf{0.193}   \\
          &           & (0.941) & (0.933) & (0.929) & (0.935) & \textbf{(0.939)} &  & (0.947) & (0.942) & (0.940) & (0.937) & \textbf{(0.941)} \\
DR        &           & 0.345   & 0.305   & 0.279   & 0.260   & \textbf{0.245}   &  & 0.264   & 0.232   & 0.211   & 0.195   & \textbf{0.183}   \\
          &           & (0.933) & (0.929) & (0.934) & (0.935) & \textbf{(0.941)} &  & (0.938) & (0.940) & (0.935) & (0.940) & \textbf{(0.942)} \\
CCT       &           & 0.334   & 0.305   & 0.282   & 0.264   & 0.249   &  & \textbf{0.263}   & 0.235   & 0.214   & 0.198   & 0.186   \\
          &           & (0.936) & (0.931) & (0.929) & (0.931) & (0.932) &  & \textbf{(0.940)} & (0.934) & (0.934) & (0.935) & (0.929) \\
  \specialrule{0pt}{0.3em}{0.em}        &           & \multicolumn{11}{c}{Setting: $b = 1.2h$}                                                             \\
\specialrule{0pt}{0.3em}{0.em} TB        &           & 0.248   & 0.225   & 0.208   & 0.194   & 0.183   &  & 0.191   & 0.169   & 0.155   & 0.144   & 0.135   \\
          &           & (0.874) & (0.866) & (0.870) & (0.856) & (0.852) &  & (0.857) & (0.865) & (0.848) & (0.855) & (0.853) \\
DB        &           & 0.220   & 0.202   & 0.188   & 0.177   & 0.167   &  & 0.176   & 0.158   & 0.144   & 0.134   & 0.126   \\
          &           & (0.848) & (0.845) & (0.844) & (0.830) & (0.824) &  & (0.838) & (0.840) & (0.833) & (0.830) & (0.834) \\
TR        &           & 0.339   & 0.299   & \textbf{0.273}   & 0.254   & 0.240   &  & 0.254   & 0.224   & 0.204   & \textbf{0.188}   & 0.176   \\
          &           & (0.940) & (0.933) & \textbf{(0.940)} & (0.939) & (0.936) &  & (0.938) & (0.941) & (0.941) & \textbf{(0.944)} & (0.934) \\
DR        &           & 0.350   & 0.305   & 0.277   & \textbf{0.258}   & 0.242   &  & 0.260   & 0.227   & 0.206   & \textbf{0.190}   & 0.178   \\
          &           & (0.941) & (0.934) & (0.942) & \textbf{(0.943)} & (0.937) &  & (0.940) & (0.941) & (0.943) & \textbf{(0.945)} & (0.941) \\
CCT       &           & 0.302   & 0.275   & 0.255   & 0.239   & 0.225   &  & 0.238   & 0.212   & 0.194   & 0.179   & 0.168   \\
          &           & (0.932) & (0.929) & (0.932) & (0.929) & (0.924) &  & (0.935) & (0.938) & (0.935) & (0.929) & (0.928) \\
 \specialrule{0pt}{0.3em}{0.em}         &           & \multicolumn{11}{c}{Setting: $b =1.5 h$}                                                             \\
\specialrule{0pt}{0.3em}{0.em} TB        &           & 0.239   & 0.218   & 0.202   & 0.189   & 0.178   &  & 0.186   & 0.166   & 0.152   & 0.141   & 0.132   \\
          &           & (0.891) & (0.888) & (0.878) & (0.864) & (0.831) &  & (0.883) & (0.875) & (0.874) & (0.877) & (0.859) \\
DB        &           & 0.221   & 0.203   & 0.189   & 0.177   & 0.168   &  & 0.176   & 0.158   & 0.145   & 0.134   & 0.126   \\
          &           & (0.869) & (0.866) & (0.854) & (0.843) & (0.814) &  & (0.857) & (0.856) & (0.853) & (0.854) & (0.849) \\
TR        &           & 0.305   & 0.273   & \textbf{0.251}   & 0.234   & 0.221   &  & 0.231   & 0.205   & 0.187   & 0.173   & 0.162   \\
          &           & (0.939) & (0.938) & \textbf{(0.942)} & (0.935) & (0.927) &  & (0.944) & (0.947) & (0.939) & (0.934) & (0.920) \\
DR        &           & 0.333   & 0.292   & 0.266   & \textbf{0.247}   & 0.232   &  & 0.247   & 0.217   & 0.197   & \textbf{0.182}   & 0.170   \\
          &           & (0.942) & (0.944) & (0.947) & \textbf{(0.948)} & (0.937) &  & (0.944) & (0.946) & (0.946) & \textbf{(0.946)} & (0.934) \\
CCT       &           & 0.276   & 0.252   & 0.233   & 0.219   & 0.207   &  & 0.217   & 0.194   & 0.177   & 0.164   & 0.154   \\
          &           & (0.932) & (0.930) & (0.933) & (0.921) & (0.904) &  & (0.939) & (0.939) & (0.929) & (0.923) & (0.911) \\ \hline
\end{tabular}
		}	
	\end{center}
\end{table}

\begin{table}[ht]
	\caption{\label{CI.case3.time}  \textcolor{black}{Average  CPU time of Model 3 with $b = 1.2h$ over 1000 simulation runs. } }
	\begin{center}
	\vspace{-0.5cm}
		\resizebox{5.5in}{!}{
       
\begin{tabular}{ccccccccccccc}
\hline
                     &                      & \multicolumn{5}{c}{$n=500$}                     &    & \multicolumn{5}{c}{$n=1000$}                 \\ \cline{3-7} \cline{9-13} 
$h$                  &                      & 0.15    & 0.18    & 0.21    & 0.24    & 0.27    &    & 0.12    & 0.15    & 0.18   & 0.21   & 0.24   \\ \hline
   \specialrule{0pt}{0.3em}{0.em}                  &                      & \multicolumn{11}{c}{Hypothesis test  for \(\tau_{\SRD}=0.04\)} \\
\specialrule{0pt}{0.3em}{0.em} Orig                 &                      & 0.005   & 0.005   & 0.005   & 0.005   & 0.005   &    & 0.006   & 0.006   & 0.005  & 0.006  & 0.007  \\
TB                   &                      & 0.005   & 0.004   & 0.005   & 0.006   & 0.005   &    & 0.006   & 0.006   & 0.006  & 0.007  & 0.006  \\
DB                   &                      & 0.011   & 0.014   & 0.014   & 0.016   & 0.019   &    & 0.022   & 0.029   & 0.035  & 0.043  & 0.052  \\
TR                   &                      & 0.005   & 0.006   & 0.006   & 0.005   & 0.006   &    & 0.006   & 0.006   & 0.007  & 0.007  & 0.008  \\
DR                   &                      & 0.011   & 0.012   & 0.014   & 0.016   & 0.019   &    & 0.018   & 0.024   & 0.034  & 0.043  & 0.056  \\
CCT                  &                      & 0.005   & 0.004   & 0.005   & 0.005   & 0.005   &    & 0.008   & 0.007   & 0.008  & 0.007  & 0.007  \\
\specialrule{0pt}{0.3em}{0.em} & \multicolumn{1}{l}{} & \multicolumn{11}{c}{Confidence interval construction}                                                     \\
\specialrule{0pt}{0.3em}{0.em} Orig                 &                      & 0.200   & 0.198   & 0.200   & 0.199   & 0.201   &    & 0.203   & 0.208   & 0.218  & 0.228  & 0.239  \\
TB                   &                      & 0.207   & 0.205   & 0.206   & 0.207   & 0.206   &    & 0.209   & 0.213   & 0.222  & 0.233  & 0.239  \\
DB                   &                      & 0.206   & 0.209   & 0.209   & 0.211   & 0.215   &    & 0.216   & 0.229   & 0.246  & 0.264  & 0.280  \\
TR                   &                      & 0.225   & 0.223   & 0.222   & 0.224   & 0.226   &    & 0.225   & 0.237   & 0.242  & 0.252  & 0.263  \\
DR                   &                      & 0.254   & 0.257   & 0.267   & 0.278   & 0.288   &    & 0.289   & 0.317   & 0.345  & 0.376  & 0.405  \\
CCT                  &                      & 0.005   & 0.004   & 0.005   & 0.005   & 0.005   &    & 0.008   & 0.007   & 0.008  & 0.007  & 0.007  \\ \hline
\end{tabular}

		}	
	\end{center}
\end{table}

\begin{table}[ht]
	\caption{\label{opt.case3} {\color{black}Comparison of average and standard deviation (in parentheses) of bandwidths $h$ and $b$, empirical sizes, empirical coverages and average interval lengths over 10000 simulation runs  for Model 3.
    The best performances, with empirical sizes close to 5\% and
empirical coverages close to 95\%, are in bold font.}}
	\begin{center}
	\vspace{-0.5cm}
		\resizebox{6.6in}{!}{
\begin{tabular}{ccccccclccccc}
\hline
\multirow{2}{*}{Method}              &                      & \multicolumn{5}{c}{$n = 500$}                                                                                                                                                             & \multicolumn{1}{c}{} & \multicolumn{5}{c}{$n = 1000$}                                                                                                                                                                          \\ \cline{3-7} \cline{9-13} 
                                     &                      & $h$                                                                           & $b$                                                                           & Size  & Coverage & Length & \multicolumn{1}{c}{} & $h$                                                                           & $b$                                                                           & Size  & Coverage & Length               \\ \hline
Orig                                 &                      & \begin{tabular}[c]{@{}c@{}}0.203\\      (0.039)\end{tabular}                  & -                                                                             & 0.092 & 0.908    & 0.198  &                      & \begin{tabular}[c]{@{}c@{}}0.188\\      (0.038)\end{tabular}                  & -                                                                             & 0.120 & 0.880    & 0.147                \\
\specialrule{0pt}{0.35em}{0.em}      & \multicolumn{1}{l}{} & \multicolumn{10}{c}{$(h,b) = (\hat   h_{\text{mse}}, \hat b_{\text{mse}})$}                                                                                                                                                                                                                                                                                                                         & \multicolumn{1}{l}{} \\
\specialrule{0pt}{0.25em}{0.em}   TB &                      & \multirow{5}{*}{\begin{tabular}[c]{@{}c@{}}0.203\\      (0.039)\end{tabular}} & \multirow{5}{*}{\begin{tabular}[c]{@{}c@{}}0.332\\      (0.063)\end{tabular}} & 0.113 & 0.887    & 0.203  &                      & \multirow{5}{*}{\begin{tabular}[c]{@{}c@{}}0.188\\      (0.038)\end{tabular}} & \multirow{5}{*}{\begin{tabular}[c]{@{}c@{}}0.319\\      (0.060)\end{tabular}} & 0.132 & 0.868    & 0.150                \\
DB                                   &                      &                                                                               &                                                                               & 0.126 & 0.874    & 0.193  &                      &                                                                               &                                                                               & 0.144 & 0.856    & 0.144                \\
TR                                   &                      &                                                                               &                                                                               & 0.060 & 0.940    & 0.246  &                      &                                                                               &                                                                               & 0.073 & 0.927    & 0.177                \\
DR                                   &                      &                                                                               &                                                                               & \textbf{0.051} & \textbf{0.949}    & 0.263  &                      &                                                                               &                                                                               & \textbf{0.060} & \textbf{0.941}    & 0.188                \\
CCT                                  &                      &                                                                               &                                                                               & 0.067 & 0.933    & 0.231  &                      &                                                                               &                                                                               & 0.084 & 0.916    & 0.169                \\
\specialrule{0pt}{0.35em}{0.em}      & \multicolumn{1}{l}{} & \multicolumn{10}{c}{{\color{black}$(h,b) = (\hat   h_{\text{mse}}, \hat h_{\text{mse}})$}}                                                                                                                                                                                                                                                                                                                         & \multicolumn{1}{l}{} \\
\specialrule{0pt}{0.25em}{0.em}   TB &                      & \multirow{5}{*}{\begin{tabular}[c]{@{}c@{}}0.203\\      (0.039)\end{tabular}} & \multirow{5}{*}{\begin{tabular}[c]{@{}c@{}}0.203\\      (0.039)\end{tabular}} & 0.128 & 0.872    & 0.218  &                      & \multirow{5}{*}{\begin{tabular}[c]{@{}c@{}}0.188\\      (0.038)\end{tabular}} & \multirow{5}{*}{\begin{tabular}[c]{@{}c@{}}0.188\\      (0.038)\end{tabular}} & 0.157 & 0.844    & 0.158                \\
DB                                   &                      &                                                                               &                                                                               & 0.141 & 0.859    & 0.192  &                      &                                                                               &                                                                               & 0.171 & 0.829    & 0.143                \\
TR                                   &                      &                                                                               &                                                                               & 0.041 & 0.959    & 0.301  &                      &                                                                               &                                                                               & \textbf{0.051} & \textbf{0.949}    & 0.222                \\
DR                                   &                      &                                                                               &                                                                               & \textbf{0.043}&\textbf{0.957}    & 0.282  &                      &                                                                               &                                                                               & \textbf{0.052} & \textbf{0.948}   & 0.210                \\
CCT                                  &                      &                                                                               &                                                                               & 0.043 & 0.958    & 0.286  &                      &                                                                               &                                                                               & 0.057 & 0.943    & 0.213                \\
\specialrule{0pt}{0.35em}{0.em}      & \multicolumn{1}{l}{} & \multicolumn{10}{c}{$(h,b) = (\hat   h_{\text{mse}}, 1.2\hat h_{\text{mse}})$}                                                                                                                                                                                                                                                                                                                      & \multicolumn{1}{l}{} \\
\specialrule{0pt}{0.25em}{0.em}   TB &                      & \multirow{5}{*}{\begin{tabular}[c]{@{}c@{}}0.203\\      (0.039)\end{tabular}} & \multirow{5}{*}{\begin{tabular}[c]{@{}c@{}}0.243\\      (0.047)\end{tabular}} & 0.118 & 0.882    & 0.210  &                      & \multirow{5}{*}{\begin{tabular}[c]{@{}c@{}}0.188\\      (0.038)\end{tabular}} & \multirow{5}{*}{\begin{tabular}[c]{@{}c@{}}0.226\\      (0.046)\end{tabular}} & 0.142 & 0.858    & 0.154                \\
DB                                   &                      &                                                                               &                                                                               & 0.134 & 0.866    & 0.192  &                      &                                                                               &                                                                               & 0.161 & 0.839    & 0.144                \\
TR                                   &                      &                                                                               &                                                                               & \textbf{0.045} & \textbf{0.955}    & 0.278  &                      &                                                                               &                                                                               & \textbf{0.054} & \textbf{0.946}    & 0.203                \\
DR                                   &                      &                                                                               &                                                                               & 0.042 & 0.958    & 0.283  &                      &                                                                               &                                                                               & \textbf{0.050} & \textbf{0.950}    & 0.205                \\
CCT                                  &                      &                                                                               &                                                                               & \textbf{0.049} & \textbf{0.951}    & 0.259  &                      &                                                                               &                                                                               & 0.061 & 0.940    & 0.192                \\
\specialrule{0pt}{0.35em}{0.em}      & \multicolumn{1}{l}{} & \multicolumn{10}{c}{$(h,b) = (\hat   h_{\text{mse}}, 1.5\hat h_{\text{mse}})$}                                                                                                                                                                                                                                                                                                                      & \multicolumn{1}{l}{} \\
\specialrule{0pt}{0.25em}{0.em}   TB &                      & \multirow{5}{*}{\begin{tabular}[c]{@{}c@{}}0.203\\      (0.039)\end{tabular}} & \multirow{5}{*}{\begin{tabular}[c]{@{}c@{}}0.304\\      (0.059)\end{tabular}} & 0.115 & 0.885    & 0.205  &                      & \multirow{5}{*}{\begin{tabular}[c]{@{}c@{}}0.188\\      (0.038)\end{tabular}} & \multirow{5}{*}{\begin{tabular}[c]{@{}c@{}}0.283\\      (0.057)\end{tabular}} & 0.137 & 0.863    & 0.151                \\
DB                                   &                      &                                                                               &                                                                               & 0.129 & 0.871    & 0.193  &                      &                                                                               &                                                                               & 0.154 & 0.846    & 0.144                \\
TR                                   &                      &                                                                               &                                                                               & \textbf{0.054} & \textbf{0.946}    & 0.256  &                      &                                                                               &                                                                               & 0.065 & 0.935    & 0.186                \\
DR                                   &                      &                                                                               &                                                                               & \textbf{0.046} & \textbf{0.954}    & 0.273  &                      &                                                                               &                                                                               & \textbf{0.057} & \textbf{0.943}    & 0.196                \\
CCT                                  &                      &                                                                               &                                                                               & 0.064 & 0.936    & 0.238  &                      &                                                                               &                                                                               & 0.076 & 0.924    & 0.176          \\ \hline     
\end{tabular}
		}	
	\end{center}
\end{table}

\section{Real data analysis}
\label{sec:real}

In this section, we apply our proposed TR and DR methods to evaluate the sharp RD effects in two real-world social study examples. We also implement the Orig, TB, DB and CCT methods, as introduced in Section~\ref{sec:sim}. To ensure a fair comparison, the MSE-optimal bandwidths $\hat h_{\text{mse}}$ and $\hat b_{\text{mse}}$  are used across all competing methods.

\subsection{Brazilian mayoral elections}
\label{sec:real_1}
The first dataset, analyzed in \cite{klavsnja2017}, contains municipal mayoral election vote shares for political parties in Brazil during the period 1996 to 2012. We aim to investigate the impact of a party’s victory in an election on its likelihood of winning the subsequent election within the same municipality. For a given party, the vote margin is defined as the difference between the party’s vote share and that of its strongest opponent in the same election. Then, the forcing variable $X_i$ represents the vote margin at the $t$-th election, and the outcome variable $Y_i$ denotes the vote margin in the subsequent $(t+1)$-th election. To clarify, we only consider municipalities where the incumbent party competes for reelection at $t+1$. The final dataset includes $n = 5460$ municipalities, divided into a treatment group and a control group. The treatment group consists of 3242 municipalities where the party won the $t$-election ($X_i \geq 0$), while the control group includes 2218 municipalities where the party lost at $t$ ($X_i < 0$). We thus focus on testing whether $\tau_{\SRD}=0$ at the cutoff point $x = 0$, which serves as the threshold between electoral loss and victory. For a detailed falsification analysis validating the RD design, see \cite{Cattaneo2020}.

Table~\ref{real.brazil} presents $p$-values and confidence intervals at $95\%$ confidence level for the competing methods, together with the selected values of $\hat h_{\text{mse}}$ and $\hat b_{\text{mse}}$. Several patterns are observable. First, all methods suggest a statistically significant and negative effect with $p$-values close to 0.  This agrees with \cite{klavsnja2017}, which concluded that becoming the incumbent party can lead to electoral losses in subsequent elections. 
Second, the three robust methods (TR, DR and CCT) produce wider confidence intervals, reflecting their improved ability to account for bias and variability in the bias estimation. Specifically, TR yields a confidence interval of $[-10.621, -3.290]$, while DR returns $[-10.720, -3.088]$.
The CCT method provides a comparable confidence interval $[-10.576, -3.340]$. Third, among the robust methods, CCT delivers a symmetric confidence interval centered at $\hat \tau_\SRD$ due to its construction,  whereas both TR and DR produce asymmetric and fully data-adaptive confidence intervals.
Notably, DR results in a slightly wider interval compared to CCT and TR, which may indicate improved empirical coverage, as evidenced by Table~\ref{opt.case3} in Section~\ref{sec:sim.RDD}.
{\color{black}Table~S7 in the supplementary material further reports the
numerical results with \(b=\hat h_{\mathrm{mse}}=14.342\). The conclusions are
similar. All methods continue to reject the null at the \(5\%\) level and imply
a negative effect of incumbency on the subsequent vote margin. As expected,
setting \(b=h\) leads to narrower confidence intervals than those obtained with
the larger $\hat b_{\text{mse}}$.}


\begin{table}[t]
	\caption{\label{real.brazil} Effect of winning at $t$ on vote margin at $t + 1$ for the incumbent party in Brazil.}
	\begin{center}
	\vspace{-0.5cm}
		\resizebox{6.0in}{!}{
\begin{tabular}{ccc|ccc|cc}
\hline
    Method & $p$-value & $95\%$ CI       &  Method    & $p$-value & $95\%$ CI         & $\hat h_{\text{mse}}$   & $\hat b_{\text{mse}}$   \\ 
\hline
Orig   & 0.000         & [-9.715, -3.354]  & CCT    & 0.000         & [-10.576, -3.340] & \multirow{3}{*}{14.342} & \multirow{3}{*}{27.912} \\
TB     & 0.000         & [-10.133, -3.760] & TR     & 0.000         & [-10.621, -3.290] &                         &                         \\
DB     & 0.000         & [-10.044, -3.719] & DR     & 0.000         & [-10.720, -3.088] &                         &                        \\ \hline
\end{tabular}
		}	
	\end{center}
\end{table}

\subsection{Turkey’s female educational attainment}
Our second dataset combines municipal mayoral election data from Turkey’s 1994 elections with educational attainment data from the 2000 Turkish Population Census, as analyzed in \cite{meyersson2014}. The goal is to examine the effect of Islamic political representation in the 1994 municipal elections on high school attainment for women whose education could have been influenced between 1994 and 2000. The matched dataset includes $n=2629$ municipalities. The forcing variable $X_i$ is the vote margin, defined as the vote percentage of the Islamic party minus that of its strongest secular opponent. Mayoral elections were determined by plurality, thus the municipalities with $X_i\geq 0$ elected an Islamic mayor (treatment group), while those with $X_i<0$ elected a secular mayor (control group). The outcome variable $Y_i$
  measures the educational attainment of women who were potentially in high school during the study period, calculated as the percentage of women aged 15 to 20 in 2000 who had completed high school by that year. For further details and validity checks on this RD design, see \cite{meyersson2014}.

Table~\ref{real.Islamic} reports the $p$-values and $95\%$ confidence intervals. While all methods suggest positive effects, the levels of significance vary. The Orig, TB and DB methods produce confidence intervals that exclude zero and yield $p$-values below 0.05. In contrast, the three robust methods provide intervals that include zero and are associated with $p$-values above 0.05, offering a more cautious interpretation of the effect of Islamic rule.  

{\color{black}To further assess the reliability of the findings, we consider three
additional pilot bandwidth choices,
\(b=\hat h_{\mathrm{mse}}, 1.2\hat h_{\mathrm{mse}}\), and
\(1.5\hat h_{\mathrm{mse}}\). Table~S8 in the supplementary material reports the corresponding \(p\)-values and \(95\%\) confidence intervals. Notably, the TB and DB methods produce inconsistent conclusions across different choices of \(b\). In particular, the confidence intervals based on
TB and DB include zero when \(b=\hat h_{\mathrm{mse}}\), but exclude zero when larger $b$ is used. This pattern aligns well with Theorem~\ref{th1}. As discussed in Remark~\ref{rm.1}, the limiting behavior of these conventional bias-corrected EL ratios depends on  \(h/b\), and thus their finite-sample performance can be sensitive to the bandwidth choices. 
 Tables~\ref{real.Islamic} and S8 reveal that this sensitivity can indeed result in over-rejection when  the \(\chi^2_1\)
distribution is used, with the associated confidence intervals
potentially undercovering and leading to unreliable inference.
By contrast, the robust methods, CCT, TR, and DR, deliver confidence intervals that include zero across all choices of \(b\). This provides
additional empirical support for the proposed robust EL methods, in the sense that the proposed new weights effectively incorporate the
variability arising from bias estimation and thereby ensure robust and valid
inference in finite samples.}


\begin{table}[t]
	\caption{\label{real.Islamic} Effect of Islamic rule on female high school educational attainment in Turkey.}
	\begin{center}
	\vspace{-0.5cm}
		\resizebox{6.0in}{!}{
\begin{tabular}{ccc|ccc|cc}
\hline
  Method   & $p$-value & $95\%$ CI      &  Method   & $p$-value & $95\%$ CI       & $\hat h_{\text{mse}}$   & $\hat b_{\text{mse}}$   \\ \hline
Orig   & 0.020      & [0.500, 6.062] & CCT    & 0.056     & [-0.080, 6.458] & \multirow{3}{*}{16.276} & \multirow{3}{*}{27.923} \\
TB     & 0.021     & [0.468, 6.033] & TR     & 0.051     & [-0.018, 6.581] &                         &                         \\
DB     & 0.024     & [0.409, 5.947] & DR     & 0.065     & [-0.189, 6.634] &                         &                            \\\hline                    
\end{tabular}
		}	
	\end{center}
\end{table}

\section{Discussion}
\label{sec:diss}
Our paper introduces a novel EL-based  strategy for nonparametric regression and RDD inference.  The key innovation lies in the development of fully data-adaptive robust weights that simultaneously correct for bias and account for its variability, thereby enabling an automatic and valid EL-based inference procedure that retains the Wilks-type chi-squared limiting distribution while exhibiting strong robustness properties.
Looking ahead, several important directions remain to be explored.

{\color{black}
First, a very important and related topic is to modify the
proposed robust EL procedures to attain a sharp coverage-error rate,
as discussed in Section~\ref{sec.second}. Our careful investigation in the
nonparametric regression setting suggests that the bias-correction term should be re-incorporated into the robust EL constraint, leading to the
following ``Version~2'' of the Taylor-expansion-based robust EL ratio,
\begin{equation*}
\begin{aligned}
\widetilde l_{\TR}(\theta) = - 2 \max \bigg\{\sum_{i = 1}^n \log(np_i) \,\bigg|\, p_i \ge 0, \,\sum_{i=1}^n p_i = 1,\, \sum_{i=1}^n p_i  W_{i,h,b}^{\star}(x) \big(Y_i - \widehat{r}_{3,b}(X_i)-\theta\big) = 0 \bigg\},
\end{aligned}
\end{equation*}
where the higher-order bias correction is defined as
$$
\widehat{r}_{3,b}(X_i) = (X_i - x) \widehat m_{2,b}^{(1)}(x) + \frac{1}{2}(X_i - x)^2 \widehat m_{2,b}^{(2)}(x).
$$
Notably, the first-order bias term within $\widehat{r}_{3,b}(X_i)$ is essential and cannot be omitted or substituted with $\widehat{r}_{1,b}(X_i)$. 
A difference-based adjustment using $\widehat r_{2,b}(X_i)$ may also be
possible, but its validity requires rigorous theoretical investigation.
This new EL framework can also be directly applicable to the sharp and fuzzy RDD scenarios. Our future work will therefore address several key objectives: (i) deriving the second-order asymptotic properties of the modified
procedure and developing a coverage-error-optimal bandwidth selector
tailored to the EL framework; (ii) conducting a thorough comparison with the bandwidth selection methods for robust confidence intervals introduced by \cite{calonico2020ej}; (iii) investigating whether simple Bartlett-type
corrections can further improve coverage accuracy. 

Another promising direction is to extend the proposed methodology to nonparametric quantile regression and quantile RDD; see, e.g., \cite{qu2019uniform} and \cite{xu2020inference}. For a sharp quantile RDD, under
the standard identification and continuity conditions, the quantile
treatment effect at level $\tau\in(0,1)$ is identified as
\[
\theta_{\SRD}(\tau)
=
\lim_{x\to0^+}Q_{Y(1)}(\tau\mid X=x)
-
\lim_{x\to0^-}Q_{Y(0)}(\tau\mid X=x),
\]
where $Q_Y(\tau\mid X=x)$ denotes the conditional $\tau$th quantile of
$Y$ given $X=x$.
Write
\[
\bG_i^{\star}(\theta,a;\tau)
=
\Big(W_{i,h,b}^{\star+}
\big\{I(Y_i\leq \theta+a)-\tau\big\},
W_{i,h,b}^{\star-}
\big\{I(Y_i\leq a)-\tau\big\}\Big)^{\T},
\]
where $W_{i,h,b}^{\star+}$ and $W_{i,h,b}^{\star-}$ are the
robust weights defined in Section~\ref{sec:4.1}. The Taylor-expansion-based robust empirical
log-likelihood function can then be formulated as
\begin{equation*}
\widetilde l_{\SRD,\TR}^{\star}(\theta,a;\tau)
=
-2\max\left\{
\sum_{i=1}^n\log(np_i)
\,\middle|\,
p_i\geq0,\quad
\sum_{i=1}^np_i=1,\quad
\sum_{i=1}^np_i\bG_i^{\star}(\theta,a;\tau)=0
\right\},
\end{equation*}
with the  empirical log-likelihood ratio
\[
l_{\SRD,\TR}^{\star}(\theta;\tau)
=
\inf_{a\in\mathcal A}
\widetilde l_{\SRD,\TR}^{\star}(\theta,a;\tau).
\]
The main methodological challenges in this extension include conducting uniform inference across multiple quantile levels within the proposed EL framework and developing the associated coverage-error-optimal bandwidth selection.

Given that empirical RDD applications often include additional covariates (see, e.g., \cite{cattaneo2023covariate,calonico2025}), a third important direction is to incorporate such auxiliary information effectively into the proposed robust EL framework to improve precision and account for potential covariate effects. These topics fall beyond the scope of the current paper and will be pursued elsewhere.

}

\section*{Supplementary material}
The supplementary material, which contains implementation details of the proposed EL procedures for RDD analysis, additional empirical results, and the higher-order representation, is available upon request.

 

\appendix
\section*{Appendix}

This appendix contains the proofs of all the main theoretical results, i.e. Theorems~\ref{th1}--\ref{th4}.

\section{Proof of Theorem~\ref{th1}}
\label{appendix-1}

We focus on proving part (b) of Theorem~\ref{th1}. Part (a) then follows by setting $\kappa = 0$.

\subsection{Asymptotic distribution of $l_{\TB}\{m(x)\}$}
\label{appendix-1-1}
 
Recall that 
$Z_i(x) = W_{i,h}(x)\{Y_i - 2^{-1}\widehat{m}_{2,b}^{(2)}(x)(X_i - x)^2 - m(x)\},$
and
$$U_1(x) = \frac{1}{n}\sum_{i=1}^n Z_i(x), \quad U_2(x) = \frac{1}{n}\sum_{i=1}^n \big\{Z_i(x)\big\}^2.$$
Using the Lagrange multiplier method, we have that
$$l_{\TB}\{m(x)\} = 2\sum_{i=1}^n \log\big\{1 + \lambda Z_i(x)\big\},$$
where $\lambda$ satisfies the equation
\begin{equation}
\label{proof-1-1}
\sum_{i=1}^{n} \frac{Z_i(x)}{1 + \lambda Z_i(x)} = 0.
\end{equation}

Analogous to the proofs in \cite{Owen2001book}, we structure the proof in four steps to establish the asymptotic distribution of $l_{\TB}\{m(x)\}$. Firstly, we demonstrate the asymptotic distribution of $U_1(x)$, i.e.,
\begin{equation}
\label{proof-1-2}
\sqrt{nh} U_1(x) \stackrel{\mathcal{D}}{\to} N\big(0, \sigma_{1,\kappa}^2(x)\big),
\end{equation}
and
\begin{equation}
\label{proof-1-3}
hU_2(x) = \sigma_0^2(x) \big\{1 + o_P(1)\big\},
\end{equation}
where the formulas for $\sigma_{1,\kappa}^2(x)$ and $\sigma_0^2(x)$ will be provided later. Notably, if $\kappa > 0$, then $\sigma_{1,\kappa}^2(x) \neq \sigma_0^2(x)$; otherwise, $\sigma_{1,0}^2(x)=\sigma_0^2(x).$ Next, we show that
\begin{equation}
\label{proof-1-4}
\max_{1 \le i \le n}|Z_i(x)| = o_P\big(n^{1/2}h^{-1/2}\big), \quad \lambda = O_P\big(n^{-1/2}h^{1/2}\big),
\end{equation}
and further establish that
\begin{equation}
\label{proof-1-5}
\lambda = U_2^{-1}(x)U_1(x) + o_P\left(n^{-1/2}h^{1/2}\right).
\end{equation}
Combining \eqref{proof-1-2}--\eqref{proof-1-5} and noting that $|\lambda|\max_{i \le n}|Z_i(x)| = o_P(1)$, by applying Taylor's expansion, i.e., $\log(1 + x) = x - 2^{-1}x^2 + o(x^2)$, we can write $l_{\TB}\{m(x)\}$ as
\begin{equation*}
\begin{aligned}
l_{\TB}\{m(x)\} &= 2\sum_{i=1}^n \log\big\{1 + \lambda Z_i(x)\big\} \\
&= 2\lambda \sum_{i=1}^n Z_i(x) - \lambda^2 \sum_{i=1}^n Z_i^2(x) \big\{1 + o_P(1)\big\} \\
&= \frac{\sigma_{1,\kappa}^2(x)}{\sigma_0^2(x)} \left\{\frac{\sqrt{nh}U_1(x)}{\sigma_{1,\kappa}(x)}\right\}^2 \big\{1 + o_P(1)\big\},
\end{aligned}
\end{equation*}
implying that $l_{\TB}\{m(x)\} \stackrel{\mathcal{D}}{\to} \gamma_{1,\kappa}\chi_1^2$, where $\gamma_{1,\kappa} = \sigma_{1,\kappa}^2(x)/\sigma_0^2(x)$.

We now provide the  proofs for equations \eqref{proof-1-2}--\eqref{proof-1-5}. To begin, we introduce some notation. Define $\bS$  as the matrix $(\mu_{j+k-2})_{1 \le j,k \le 3}$, and let $\mu_j = \int_{-1}^1 u^j K(u) \, du,$ $\nu_j =  \int_{-1}^1 u^jK(u)^2 \, du$ for $j = 0, \ldots, 4$. Note that, in particular, $\mu_0 = 1$ and $\mu_1 = \mu_3 = 0$. Additionally, we define $w_{i,h}(x) = \mu_2 K_h(X_i - x) $ for $i = 1, \ldots, n$.

We first demonstrate that the asymptotic distribution of 
$U_1(x).$ Recall that 
$$
\widehat m_{2,b}^{(2)}(x) =  \frac{2}{nb^2} \sum_{i=1}^n W_{i,2,2,b}(x) Y_i,
$$ where
 $ W_{i,2,2,b}(x) = \bee^{\T}_{3} \bS^{-1}_{2,b} \{1, (X_i - x)/b, (X_i-x)^2/b^2\}^\T K_{b}(X_i-x).$  Given that $m(x)$ is continuously differentiable up to the third order in the neighborhood of $x,$ it follows from Section 3.2.2 in \cite{fan1996local} that
  \begin{equation*}
\widehat m_{2,b}^{(2)}(x) - m^{(2)}(x) = O_P(b) + \frac{2}{n b^2 f(x)} \sum_{i=1}^n \widetilde w_{i,b}(x) \varepsilon_i \big\{1 + o_P(1)\big\}, 
 \end{equation*}
 where $\widetilde w_{i,b}(x) = \bee^{\T}_{3} \bS^{-1} \{1, (X_i - x)/b, (X_i-x)^2/b^2\}^\T K_{b}(X_i-x).$ This further implies that 
 $\widehat m_{2,b}^{(2)}(x) - m^{(2)}(x) = O_P\{b + (nb^5)^{-1/2}\}.$

Note that $\sum_{i=1}^n (X_i - x)W_{i,h}(x) = 0,~W_{i,h}(x) = w_{i,h}(x) f(x) \{1 + o_P(1)\},$ and
$$
\sum_{i=1}^n (X_i-x)^2 W_{i,h}(x)= nh^2f^2(x) \mu_2^2\big\{1 + o_P(1)\big\}.
$$
The term $U_1(x)$ can be reformulated as 
\begin{equation*}
\begin{aligned}
U_1(x) 
& = \frac{1}{n}\sum_{i=1}^n W_{i,h}(x)\left\{r(X_i) + \varepsilon_i\right\} - \frac{1}{2n}\sum_{i=1}^n W_{i,h}(x)(X_i-x)^2 \big\{\widehat m_{2,b}^{(2)}(x) - m^{(2)}(x)\big\}\\
& = O_P\big(h^3 + h^2b\big) + \frac{1}{n}\sum_{i=1}^n W_{i,h}(x)\varepsilon_i- \frac{h^2 \mu_2^2 f(x)}{nb^2}\sum_{i=1}^n \widetilde w_{i,b}(x) \varepsilon_i \big\{1 + o_P(1)\big\}\\
& = O_P\big(h^3 + h^2b\big) + \frac{f(x)}{n}\sum_{i=1}^n \big\{ w_{i,h}(x)- \kappa^2 \mu_2^2 \widetilde w_{i,b}(x) \big\} \varepsilon_i\big\{1 + o_P(1)\big\}.
\end{aligned}
\end{equation*}
Given that $\sqrt{nh} (h^3 + h^2 b) \to 0$ since  $nh^5b^2 \to 0$ and $h/b \to \kappa \in [0,1],$ by the classic central limit theorem, we obtain (\ref{proof-1-2}), where
$$
\sigma_{1,\kappa}^2(x) = \underset{h/b \to \kappa, b \to 0}{\lim} h E\big\{ w_{i,h}(x)- \kappa^2 \mu_2^2\widetilde w_{i,b}(x) \big\}^2\sigma^2(x) f^2(x),
$$
and this limit exists. 

We write the term $U_2(x)$ as 
\begin{equation*}
\begin{aligned}
U_{2}(x) = \frac{1}{n}\sum_{i=1}^n \big\{W_{i,h}(x)\big\}^2\big\{\varepsilon_i + o_P(1)\big\}^2 = \frac{1}{n}\sum_{i=1}^n w_{i,h}^2(x)\varepsilon_i^2 f^2(x)\big(1 + o_P(1)\big).
\end{aligned}
\end{equation*}
This implies $hU_2(x) = \sigma^2_{0}(x) \{1+o_P(1)\}$, where $
\sigma_0^2(x) = \mu_2^2\nu_0\sigma^2(x)f^3(x),
$ 
which is equal to $\sigma_{1,0}^2(x),$ and hence (\ref{proof-1-3}) follows.

Now we turn to show that $\max_{1 \le i\le n}|Z_i(x)| = o_P(n^{1/2} h^{-1/2}).$  Observe that $|S_{2,h}(x)| +  |S_{1,h}(x)| = O_P(1),$ we have 
\begin{equation*}
\begin{aligned}
  \max_{1 \le i\le n}|Z_i(x)| \le & \max_{1 \le i \le n}\big\{K_h(X_i-x)\{|\varepsilon_i| + o(1)\}\big\}\cdot O_P(1)\\
  & + \max_{1 \le i \le n}\big\{ K_h(X_i-x)\big\}  O_P\big\{ h^2b + h^2(nb^5)^{-1/2}\big\}\\
  \le & \max_{1 \le i \le n}\big\{K_h(X_i-x)\{|\varepsilon_i| + O(1)\}\big\}\cdot O_P(1).
\end{aligned}
\end{equation*}
Since $E \{K_h^4(X_i-x)\varepsilon_i^4\} = O(h^{-3}),$ it follows from Markov inequality that 
$$
\max_{1 \le i \le n}|K_h(X_i-x)\varepsilon_i|  = o_P\big(n^{1/2} h^{-1/2}\big).
$$
which further implies that 
$\max_{1 \le i\le n}|Z_i(x)| = o_P(n^{1/2} h^{-1/2}).$

Next, we show that $\lambda = O_P(n^{-1/2}h^{1/2}).$ 
From Equation (\ref{proof-1-1}), we have
\begin{equation*}
\sum_{i=1}^{n} Z_i(x)=
\sum_{i=1}^{n} \frac{\lambda Z_i^2(x)}{1+\lambda Z_i(x)}.
\end{equation*}
Since each $p_i \ge 0,$ $1+\lambda Z_i(x) >0,$ and therefore, 
\begin{equation*}
\frac{1}{n}\sum_{i=1}^{n} \frac{ Z_i^2(x)}{1+\lambda Z_i(x)} \ge \frac{U_2(x)}{1+ |\lambda| \max_{i\le n}|Z_i(x)|}.
\end{equation*}
This inequality implies that 
\begin{equation*}
|\lambda| U_2(x) \le \big(1+ |\lambda| \max_{i\le n}|Z_i(x)|\big) \bigg|\frac{1}{n} \sum_{i=1}^{n} Z_i(x)\bigg|.
\end{equation*}
Thus we conclude that
\begin{equation*}
\begin{aligned}
|\lambda| \sigma^2_2(x) \{1+o_P(1)\} & \le h\cdot O_P\big\{(nh)^{-1/2}\big\} + |\lambda| \cdot h \cdot o_P\big(n^{1/2}h^{-1/2}\big)\cdot O_P\big\{(nh)^{-1/2}\big\}\\
& \le h\cdot O_P\big\{(nh)^{-1/2}\big\} + |\lambda| \cdot o_P(1),
\end{aligned}
\end{equation*}
which in turn implies that 
$|\lambda| = O_P(n^{-1/2}h^{1/2}).$

Finally, we show that $\lambda = U_2^{-1}U_1(x) + o_P(n^{-1/2}h^{1/2}).$
First, we observe that
\begin{equation*}
\begin{aligned}
0 & = \sum_{i=1}^{n} \frac{Z_i(x)}{1+\lambda Z_i(x)}
  & =  \sum_{i=1}^{n} Z_i(x) - \lambda \sum_{i=1}^{n} Z_i^2(x) + \sum_{i=1}^{n} \frac{\lambda^2 Z_i^3(x)}{1+\lambda Z_i(x)}.
\end{aligned}
\end{equation*}
Note that $|\lambda| \max_{i\le n}|Z_i(x)| = o_P(1).$ Hence it follows that 
\begin{equation*}
\left|\sum_{i=1}^{n} \frac{\lambda^2 Z_i^3(x)}{1+\lambda Z_i(x)} \right| \le o_P(|\lambda|)\cdot\sum_{i=1}^{n} \frac{Z_i^2(x)}{1-o_P(1)} 
\le o_P(|\lambda|) \sum_{i=1}^{n} Z_i^2(x) \big\{1 + o_P(1)\big\}.
\end{equation*}
As a result, we obtain that
\begin{equation*}
\begin{aligned}
  \lambda  =  \bigg\{\sum_{i=1}^{n} Z_i^2(x) \bigg\}^{-1}\sum_{i=1}^{n} Z_i(x) + o_P(|\lambda|) =  U_2^{-1}U_1(x) + o_P\big(n^{-1/2}h^{1/2}\big).
\end{aligned}
\end{equation*}
The proof for $l_{\TB}\{m(x)\}$ is complete. 

\subsection{Asymptotic distribution of $l_{\DB}\{m(x)\}$}
 \label{appendix-1-2}
 Define
$
\widetilde Z_i(x) = W_{i,h}(x)\{Y_i- m(x) - (\widehat m_{2,b}(X_i) - \widehat m_{2,b}(x))\},
$ and let
$$
\widetilde U_1(x) = \frac{1}{n}\sum_{i=1}^n \widetilde Z_i(x), \quad \widetilde U_2(x) = \frac{1}{n}\sum_{i=1}^n \big\{\widetilde Z_i(x)\big\}^2.$$
Through the Lagrange multiplier approach, we have that
$$l_{\DB}\{m(x)\} = 2\sum_{i=1}^n \log\big\{1 + \lambda \widetilde Z_i(x)\big\},$$
where $\lambda$ satisfies the equation
\begin{equation*}
\sum_{i=1}^{n} \frac{\widetilde Z_i(x)}{1 + \lambda \widetilde Z_i(x)} = 0.
\end{equation*}

Similar to the proof techniques outlined in Section \ref{appendix-1-1}, it can be straightforwardly shown that
$$
h \widetilde{U}_2(x) = \sigma_0^2(x) \big\{1 + o_P(1)\big\}, \quad \max_{1 \le i \le n} |\widetilde{Z}_i(x)| = o_P\big(n^{1/2}h^{-1/2}\big), \quad \lambda = O_P\left( n^{-1/2}h^{1/2} \right),
$$
and furthermore, it can be established that
$
\lambda = \widetilde{U}_2^{-1}(x) \widetilde{U}_1(x) + o_P( n^{-1/2}h^{1/2} ).
$
These derivations are therefore omitted here for brevity. If we establish that
\begin{equation}
\label{proof-2-1}
\sqrt{nh} \widetilde U_1(x) \stackrel{\mathcal{D}}{\to} N\big(0, \sigma_{2,\kappa}^2(x)\big),
\end{equation}
then we can write $l_{\DB}\{m(x)\}$ as follows:
\begin{equation*}
\begin{aligned}
l_{\DB}\{m(x)\} 
&= 2\lambda \sum_{i=1}^n \widetilde Z_i(x) - \lambda^2 \sum_{i=1}^n \big\{\widetilde Z_i(x)\big\}^2 \big\{1 + o_P(1)\big\} \\
&= \frac{\sigma_{2,\kappa}^2(x)}{\sigma_0^2(x)} \left\{\frac{\sqrt{nh} \widetilde U_1(x)}{\sigma_{2,\kappa}(x)}\right\}^2 \big\{1 + o_P(1)\big\},
\end{aligned}
\end{equation*}
implying that $l_{\DB}\{m(x)\} \stackrel{\mathcal{D}}{\to} \gamma_{2,\kappa}\chi_1^2$, where $\gamma_{2,\kappa} = \sigma_{2,\kappa}^2(x)/\sigma_0^2(x)$.
Therefore, in the following discussion, we will focus on deriving the asymptotic distribution of $\widetilde{U}_1(x)$, which corresponds to (\ref{proof-2-1}).

 First, it follows from Lemma 1 of \cite{fan1999} that 
 \begin{equation*}
 \begin{aligned}
\widehat m_{2,b}(x) - m(x)  =  \frac{1}{nf(x)} \sum_{i=1}^n \breve w_{i,b}(x) \varepsilon_i
 + O_P\left( b^3 +  \frac{b^2 \log^{1/2} n}{n^{1/2}b^{1/2}} + \frac{\log n}{nb}\right), 
\end{aligned}
 \end{equation*}
 holds uniformly over $ x \in [0,1]$, where $\breve w_{i,b}(x) = \bee^{\T}_{1} \bS^{-1} \{1, (X_i - x)/b, (X_i-x)^2/b^2\}^\T K_{b}(X_i-x).$
Given that $b^2 \log n \to 0$ and $nh^3b^4 \to 0,$ for each $X_i$ and $x,$ 
\begin{equation*}
\widehat m_{2,b}(X_i) - \widehat m_{2,b}(x) = m(X_i) - m(x)  
+ \frac{1}{nf(x)} \sum_{k=1}^n \big\{\breve w_{k,b}(X_i) - \breve w_{k,b}(x)\big\} \varepsilon_k
 + o_P\left\{ (nh)^{-1/2}\right\},
\end{equation*}
The term $\widetilde U_1(x)$ can then be written as follows:
\begin{equation*}
\begin{aligned}
\widetilde U_1(x) 
& = \frac{f(x)}{n}\sum_{i=1}^n w_{i,h}(x)\varepsilon_i  - \frac{1}{n} \sum_{i=1}^n \frac{1}{n} \sum_{k=1}^n w_{k,h}(x)\big\{ \breve w_{i,b}(X_k) - \breve w_{i,b}(x)\} \varepsilon_i
+ o_P\left\{ (nh)^{-1/2}\right\}.
\end{aligned}
\end{equation*}
Denote an equivalent kernel $\breve K(t)$ by
$$
\breve K(t) = \bee^{\T}_{1} \bS^{-1} \int_{-1}^1  \{1,t - \kappa u,(t - \kappa u)^2 \}^\T K(u) K(t - \kappa u) {\rm d}u,
$$
and $\breve K_{i,b}(x) = \breve K\{(X_i - x)/b\}/b,$ for $i =1,\ldots,n.$  We can show that for each $X_i$ and $x,$
\begin{equation*}
\begin{aligned}
\frac{1}{n}\sum_{k=1}^n E_{X_k} \big\{w_{k,h}(x) \breve w_{i,b}(X_k)\big\} = \mu_2\breve K_{i,b}(x) f(x) \{1 + o(1)\},
\end{aligned}
\end{equation*}
Moreover, for each $X_i$ and $x$, the corresponding variance is
    $$
    \frac{1}{n^2}\sum_{j=1}^n \var_{X_j} \big\{w_{j,h}(x) \breve w_{i,b}(X_k)\big\} = O\big\{(nhb^2)^{-1}\big\} = o\big(h^{-1}\big).
    $$
since $nb^2 \to \infty$ as $n \to \infty.$ Note that $n^{-1}\sum_{k=1}^n w_{k,h}(x) = \mu_2f(x) \{1 + o_P(1)\}.$ Therefore, 
\begin{equation*}
\begin{aligned}
\frac{1}{n}\sum_{k=1}^n w_{k,h}(x)\big\{ \breve w_{i,b}(X_k) - \breve w_{i,b}(x) \big\}  = \mu_2f(x) \big\{\breve K_{i,b}(x) - \breve w_{i,b}(x)\big\}\big\{1 + o_P(1)\big\}.
\end{aligned}
\end{equation*}
This yields that 
\begin{equation*}
\begin{aligned}
\widetilde U_1(x) = \frac{f(x)}{n}\sum_{i=1}^n \big\{w_{i,h}(x) - \mu_2 \breve K_{i,b}(x) + \mu_2\breve w_{i,b}(x) \big\}\varepsilon_i + o_P\big\{(nh)^{-1/2}\big\}.
\end{aligned}
\end{equation*}
Then, by the classic central limit theorem, we obtain that
\begin{equation}
\nonumber
\sqrt{nh} \widetilde U_1(x) \stackrel{\mathcal{D}}{\longrightarrow} N\big(0, \sigma_{2,\kappa}^2(x)\big),
\end{equation}
where $\sigma_{2,\kappa}^2(x) = \underset{h/b \to \kappa, b \to 0}{\lim} h E\big\{ w_{i,h}(x) - \mu_2 \breve K_{i,b}(x) + \mu_2\breve w_{i,b}(x) \big\}^2\sigma^2(x)f^2(x),$ and this limit exists.   
It is obvious that when $\kappa =0,$ $\sigma_{2,0}^2(x) = \sigma^{2}_0(x).$


The proof for $l_{\DB}\{m(x)\}$ is complete.  $\hfill\Box$

\section{Proof of Theorem~\ref{th2}}

\subsection{Asymptotic distribution of $l_{\TR}\{m(x)\}$}
\label{appendix-2-1}

Define $Z_i^\star(x) = W_{i,h,b}^\star(x) \big\{Y_i - m(x)\big\},$ and let
$$
U_1^\star(x) = \frac{1}{n}\sum_{i=1}^n Z_i^\star (x), \quad U_2^\star(x) = \frac{1}{n}\sum_{i=1}^n \big\{Z_i^\star(x)\big\}^2.$$
By the Lagrange multiplier method, we have that
$$l_{\TR}\{m(x)\} = 2\sum_{i=1}^n \log\big\{1 + \lambda Z_i^\star(x)\big\},$$
where $\lambda$ is determined by the equation
\begin{equation*}
\sum_{i=1}^{n} \frac{ Z_i^\star(x)}{1 + \lambda  Z_i^\star(x)} = 0.
\end{equation*}
Following similar proof techniques as in Section \ref{appendix-1-1}, it can be readily shown that
$$
\quad \max_{1 \le i \le n} |Z_i^\star(x)| = o_P\big(n^{1/2}h^{-1/2}\big), \quad \lambda = O_P\left( n^{-1/2}h^{1/2} \right),
$$
Since $ U_1^\star(x) $ is identical to $U_1(x)$ in Section \ref{appendix-1-1}, we have $\sqrt{nh} U_1^\star(x) \stackrel{\mathcal{D}}{\to} N\big(0, \sigma_{1,\kappa}^2(x)\big).$ Therefore, these derivations are omitted here. 
Suppose we can demonstrate that 
\begin{equation}
\label{proof-th2-1}
h U_2^\star(x) = \sigma_{1,\kappa}^2(x) \big\{1 + o_P(1)\big\}.
\end{equation}
It immediately follows that
$
\lambda = \left\{U_2^{\star}(x)\right\}^{-1} U_1^\star(x) + o_P\left( n^{-1/2}h^{1/2} \right).
$
We can then write $l_{\TR}\{m(x)\}$ as follows:
\begin{equation*}
\begin{aligned}
l_{\TR}\{m(x)\} 
&= 2\lambda \sum_{i=1}^n  Z_i^{\star}(x) - \lambda^2 \sum_{i=1}^n \big\{ Z_i^{\star}(x)\big\}^2 \big\{1 + o_P(1)\big\} \\
&= \left\{\frac{\sqrt{nh} U_1^{\star}(x)}{\sigma_{1,\kappa}(x)}\right\}^2 \big\{1 + o_P(1)\big\}\stackrel{\mathcal{D}}{\to} \chi_1^2.
\end{aligned}
\end{equation*}
We thus focus on deriving (\ref{proof-th2-1}) in the remainder of the proof.

Write $r_i(x) = m(X_i) - m(x)$. Then
\begin{equation*}
\begin{aligned}
Z_i^\star(x)  =~ & W_{i,h,b}^\star(x)  \big\{Y_i - m(x)\big\}\\ 
 =~ & W_{i,h}(x) \big\{r_i(x) + \varepsilon_i\big\} - \frac{1}{n}\sum_{k=1}^n W_{k,h}(x)(X_k-x)^2 b^{-2}W_{i,2,2,b}(x)\big\{r_i(x)+ \varepsilon_i\big\} \\
 =~ & \Big\{W_{i,h}(x) \varepsilon_i  - \frac{1}{n}\sum_{k=1}^n W_{k,h}(x)(X_k-x)^2 b^{-2}W_{i,2,2,b}(x) \varepsilon_i \Big\}\\
   & + W_{i,h}(x)r_i(x) - \frac{1}{n}\sum_{k=1}^n W_{k,h}(x)(X_k-x)^2 b^{-2}W_{i,2,2,b}(x)r_i(x)\\
\equiv~ & Z_{i,1}(x) + Z_{i,2}(x) + Z_{i,3}(x).
\end{aligned}
\end{equation*}
Let $U_{21}^\star(x) = n^{-1}\sum_{i=1}^n \{Z_{i,1}(x)\}^2.$ It is simple to show that 
\begin{equation*}
\begin{aligned}
hU_{21}^\star(x)
= \frac{h}{n}\sum_{i=1}^n \left\{w_{i,h}(x)  - \kappa^2 \mu_2^2 \widetilde w_{i,b}(x) \right\}^2 \varepsilon_i^2 f^2(x) \big\{1 + o_P(1)\big\}
= \sigma^2_{1,\kappa}(x) \big\{1 + o_P(1)\big\}.
\end{aligned}
\end{equation*}
Note that $r_i(x) = m(X_i) - m(x) = O(|X_i-x|)$ around $x$. We have 
\begin{equation*}
\begin{aligned}
hU_{22}^\star(x) =~ &  \frac{h}{n}\sum_{i=1}^n \big\{Z_{i,2}(x)\big\}^2  
= \frac{h}{n}\sum_{i=1}^n \big\{w_{i,h}(x)\big\}^2 r_i^{2}(X_i) f^2(x)\big\{1 + o_P(1)\big\} = O_P\big(h^2\big),\\
hU_{23}^\star(x) =~ &  \frac{h}{n}\sum_{i=1}^n \big\{Z_{i,3}(x)\big\}^2  
= \frac{h}{n}\sum_{i=1}^n \kappa^4 c_1^2 \big\{\widetilde w_{i,b}(x)\big\}^2 r_i^{2}(X_i) f^2(x)\big\{1 + o_P(1)\big\} = O_P\big(b^2\big).
\end{aligned}
\end{equation*}
Therefore, by the Cauchy--Schwarz inequality, we show that 
\begin{equation*}
\begin{aligned}
&\frac{h}{n}\bigg|\sum_{i=1}^n Z_{i,1}(x) Z_{i,2}(x)\bigg| \le \big\{hU_{21}^\star(x)\big\}^{1/2}\big\{hU_{22}^\star(x)\big\}^{1/2} = O_P(h),\\
&\frac{h}{n}\bigg|\sum_{i=1}^n Z_{i,1}(x) Z_{i,3}(x)\bigg| \le \big\{hU_{21}^\star(x)\big\}^{1/2}\big\{hU_{23}^\star(x)\big\}^{1/2} = O_P(b),\\
&\frac{h}{n}\bigg|\sum_{i=1}^n Z_{i,2}(x) Z_{i,3}(x)\bigg| \le \big\{hU_{22}^\star(x)\big\}^{1/2}\big\{hU_{23}^\star(x)\big\}^{1/2} = O_P(hb), 
\end{aligned}
\end{equation*}
and hence, 
\begin{equation*}
hU_{2}^\star(x) = \sigma^2_{1,\kappa}(x) \big\{1 + o_P(1)\big\}.
\end{equation*}
The proof for $l_{\TR}\{m(x)\}$ is complete. 

\subsection{Asymptotic distribution of $l_{\DR}\{m(x)\}$}

Let $Z_i^\diamond(x) = W_{i,h,b}^{\diamond}(x) \left\{Y_i - m(x)\right\}, i=1,\ldots,n$ and
$
U_{2}^\diamond(x) = n^{-1} \sum_{i=1}^n \left\{Z_i^\diamond(x)\right\}^2. 
$
Following similar arguments as in Section~\ref{appendix-2-1}, it suffices to show that
$$
hU_{2}^\diamond(x) = \sigma^2_{2,\kappa}(x) \big\{1 + o_P(1)\big\}. 
$$

Observe that 
\begin{equation*}
\begin{aligned}
Z_i^\diamond(x)  
 =~ & W_{i,h}(x) \big\{r_i(x) + \varepsilon_i\big\} - \frac{1}{n}\sum_{k=1}^n W_{k,h}(x) \big\{W_{i,0,2,b}(X_k) - W_{i,0,2,b}(x)\big\}\big\{r_i(x)+ \varepsilon_i\big\} \\
 =~ & \Big[W_{i,h}(x) \varepsilon_i  - \frac{1}{n}\sum_{k=1}^n W_{k,h}(x) \big\{W_{i,0,2,b}(X_k) - W_{i,0,2,b}(x)\big\} \varepsilon_i \Big]\\
   & + W_{i,h}(x)r_i(x) - \frac{1}{n}\sum_{k=1}^n W_{k,h}(x) \big\{W_{i,0,2,b}(X_k) - W_{i,0,2,b}(x)\big\}r_i(x)\\
\equiv~ & Z_{i,1}^\diamond(x) +  Z_{i,2}^\diamond(x) +  Z_{i,3}^\diamond(x).
\end{aligned}
\end{equation*}
with $r_i(x) = m(X_i) - m(x)$. Note that $W_{i,h} = w_{i,h}f(x)\{1 + o_P(1)\},$ and
\begin{equation*}
\begin{aligned}
\frac{1}{n}\sum_{k=1}^n W_{k,h}(x) \big\{W_{i,0,2,b}(X_k) - W_{i,0,2,b}(x)\big\}  = \big\{\breve K_{i,b}(x) - \breve w_{i,b}(x)\big\} \mu_2 f(x)\big\{1 +o_P(1)\big\},
\end{aligned}
\end{equation*}
where $\breve K_{i,b}(x)$ is defined in Section \ref{appendix-1-2}.
Then,
\begin{equation*}
\begin{aligned}
hU_{21}^\diamond(x) =~ &  \frac{h}{n}\sum_{i=1}^n \big\{Z_{i,1}^\diamond(x) \big\}^2 \\  
=~ &  \frac{h}{n}\sum_{i=1}^n \big\{w_{i,h}(x) - \mu_2\breve K_{i,b}(x) + \mu_2 \breve w_{i,b}(x) \big\}^2 \varepsilon_i^2 f^2(x) \big\{1 + o_P(1)\big\}\\
=~ &  \sigma^2_{2,\kappa}(x) \{1 + o_P(1)\}.
\end{aligned}
\end{equation*}
Since $r_i(x) = m(X_i) - m(x) = O(|X_i-x|)$ around $x$, we have 
\begin{equation*}
\begin{aligned}
h U_{22}^\diamond(x) = & \frac{h}{n}\sum_{i=1}^n \big\{Z_{i,2}^\diamond(x)\big\}^2 = \frac{h}{n}\sum_{i=1}^n \big\{w_{i,h}(x)\big\}^2 r_i^{2}(X_i) f^2(x)\big\{1 + o_P(1)\big\} = O_P\big(h^2\big),\\
hU_{23}^\diamond(x) =~ &  \frac{h}{n}\sum_{i=1}^n \big\{Z_{i,3}^\diamond(x)\big\}^2 = \frac{h}{n}\sum_{i=1}^n  \big\{ \breve K_{i,b}(x) - \breve w_{i,b}(x)\big\}^2 r_i^{2}(X_i) \mu_2^2f^2(x)\big\{1 + o_P(1)\big\} = O_P\big(b^2\big).
\end{aligned}
\end{equation*}
Moreover, by Cauchy--Schwartz inequality, we obtain that 
\begin{eqnarray*}
\frac{h}{n}\sum_{i=1} \left\{Z_{i,1}^\diamond(x)Z_{i,2}^\diamond(x) + Z_{i,1}^\diamond(x)Z_{i,3}^\diamond(x) + Z_{i,2}^\diamond(x)Z_{i,3}^\diamond(x)\right\} = O_P(h + b).
\end{eqnarray*}
Hence, it follows that
\begin{equation*}
h U_{2}^\diamond(x) = \sigma^2_{2,\kappa}(x) \big\{1 + o_P(1)\big\}.
\end{equation*}
The proof for $l_{\DR}\{m(x)\}$ is complete. $\hfill\Box$

\section{Proofs of Theorems~\ref{th3} and \ref{th4}}

We begin by introducing some notation and useful asymptotic results. Recall
$$
\bZ_i^{\star}(\theta, a) = \big(W_{i,h,b}^{\star+}(Y_i - \theta - a), W_{i,h,b}^{\star-}(Y_i - a)\big)^\T.
$$
Define $\bW_i = (W_{i,h,b}^{\star+}, W_{i,h,b}^{\star-})^\T$.
Let $\widehat{\bW}_n = n^{-1} \sum_{i=1}^n \bW_i$, $\bU_n = n^{-1} \sum_{i=1}^n \bZ_i^{\star}(\tau_{\SRD}, \mu_-)$, and
$$
\widehat{\bV}_n = \frac{1}{n} \sum_{i=1}^n \bZ_i^{\star}(\tau_{\SRD}, \mu_-)\bZ_i^{\star}(\tau_{\SRD}, \mu_-)^\T.
$$
Note that the first and second components of $\bU_n$ correspond to the bias-corrected local linear estimators $\widehat{\mu}_+$ and $\widehat{\mu}_-$, respectively, as defined in \cite{Calonico2014}. Let $c_k = \int_0^1 K(u) u^k {\rm d}u$ for $k = 0, 1, 2$. By Lemma A1 and Theorem A1 in \cite{Calonico2014}, under Conditions \ref{cond3}--\ref{cond_sharp-rdd-2}, it follows that
\begin{equation}
\label{proof-3-0}
\begin{aligned}
\sqrt{nh} \bU_n \stackrel{\mathcal{D}}{\to} N(0, \bV), \quad h \widehat{\bV}_n = \bV \big\{1 + o_P(1)\big\},
\end{aligned}
\end{equation}
where $\bV = (c_2 c_0 - c_1^2)^2 f(0)^2\diag(V_+, V_-)$, and the formulas for $V_+$ and $V_-$ are given by $\mathcal{V}^{bc}_{+,0,1,2}(h,b)$ and $\mathcal{V}^{bc}_{-,0,1,2}(h,b)$ in Theorem A1(V) of \cite{Calonico2014}, respectively. Furthermore, it also holds that
$
\widehat{\bW}_n = \bW \{1 + o_P(1)\}$ with $\bW = (c_2 c_0 - c_1^2) f(0) (1,1)^\T
$.

Now we proceed to prove that $l_{\SRD,\TR}(\tau_{\SRD}) \stackrel{\mathcal{D}}{\to } \chi^2_1.$ Let $\hat{\mu}_- = \arg \min_{a \in \mathcal{A}} l_{\SRD,\TR}(\tau_{\SRD}, a)$. The ratio $l_{\SRD,\TR}(\tau_{\SRD})$ is defined as
$$
l_{\SRD,\TR}(\tau_{\SRD}) = 2 \sum_{i=1}^n \log \big\{ 1 + \widehat \blambda^\T \bZ_i^{\star}(\tau_{\SRD}, \hat{\mu}_-) \big\},
$$
where $\widehat \blambda \in \mathcal{R}^2$, and $( \widehat \blambda^\T, \hat{\mu}_-)^\T$ is the solution to the following equations:
\begin{align}
\label{proof-3-1}
\sum_{i=1}^n \frac{\bZ_i^{\star}(\tau_{\SRD}, a)}{1 + \blambda^\T \bZ_i^{\star}(\tau_{\SRD}, a)} &= 0, \\
\label{proof-3-2}
\sum_{i=1}^n \frac{\blambda^\T \bW_i}{1 + \blambda^\T \bZ_i^{\star}(\tau_{\SRD}, a)} &= 0.
\end{align}
Suppose initially that $\hat{\mu}_- - \mu_- = O_P\{(nh)^{-1/2}\}$. We will verify this assertion at the end of this section. Together with Equation (\ref{proof-3-1}), we obtain
$$
\max_{1 \le i \le n} \|\bZ_i^{\star}(\tau_{\SRD}, \hat{\mu}_-)\|_2 = \max_{1 \le i \le n} \|\bZ_i^{\star}(\tau_{\SRD}, \mu_-)\|_2 + \max_{1 \le i \le n} \|\bW_i\|_2 |\widehat{\mu}_- - \mu_-| = o_P\big(n^{1/2} h^{-1/2}\big),
$$
and, consequently,
$
\|\widehat \blambda\|_2 = O_P(n^{-1/2} h^{1/2}),
$
where $\|\cdot\|_2$ denotes the Euclidean norm. By performing a Taylor expansion for Equations (\ref{proof-3-1}) and (\ref{proof-3-2}) around $(0, 0, \mu_-)^\T$, we obtain
\begin{equation}
\label{proof-3-3}
\begin{aligned}
\widehat{\bW}_n (\hat{\mu}_- - \mu_-) + \widehat{\bV}_n \widehat{\blambda} &= \bU_n + o_P\big(n^{-1/2} h^{-1/2}\big),\\
\widehat{\bW}_n^\T \widehat{\blambda} &= o_P\big(n^{-1/2} h^{1/2}\big).
\end{aligned}
\end{equation}
Let $\bP_n = \bI_2 - \widehat{\bV}_n^{-1/2} \widehat{\bW}_n \left( \widehat{\bW}_n^\T \widehat{\bV}_n^{-1} \widehat{\bW}_n \right)^{-1} \widehat{\bW}_n^\T \widehat{\bV}_n^{-1/2},$  where $\bI_2$ is the $2 \times 2$ identity matrix. Then, it follows from (\ref{proof-3-3}) that
\begin{equation}
\label{proof-3-6}
\begin{aligned}
\widehat{\bV}_n^{-1/2} \bW_n (\hat{\mu}_- - \mu_-) &= (\bI_2 - \bP_n) \widehat{\bV}_n^{-1/2} \bU_n + o_P\big(n^{-1/2}\big)\\
\widehat{\blambda} &= \widehat{\bV}_n^{-1/2} \bP_n \widehat{\bV}_n^{-1/2} \bU_n + o_P\big(n^{-1/2} h^{1/2}\big).
\end{aligned}
\end{equation}

Note that $n^{-1} \sum_{i=1}^n \bZ_i^{\star}(\tau_{\SRD}, \widehat{\mu}_-) = \bU_n - \widehat{\bW}_n (\widehat{\mu}_- - \mu_-)$, $h \widehat{\bV}_n = \bV \{1 + o_P(1)\}$, $\bP_n = \bP \{1 + o_P(1)\}$ with
$
\bP = \bI_2 - \bV^{-1/2} \bW \left( \bW^\T \bV^{-1} \bW \right)^{-1} \bW^\T \bV^{-1/2}.
$
Combining these results with (\ref{proof-3-6}), $l_{\SRD,\TR}(\tau_{\SRD})$ can be re-expressed as
\begin{eqnarray*}
l_{\SRD,\TR}(\tau_{\SRD}) &=& 2 \sum_{i=1}^n \log \big\{ 1 + \widehat{\blambda}^\T \bZ_i^{\star}(\tau_{\SRD}, \hat{\mu}_- ) \big\} \\
&=& 2 \widehat{\blambda}^\T \sum_{i=1}^n \bZ_i^{\star}(\tau_{\SRD}, \hat{\mu}_- ) - \widehat{\blambda}^\T \sum_{i=1}^n \bZ_i^{\star}(\tau_{\SRD}, \hat{\mu}_- ) \bZ_i^{\star}(\tau_{\SRD}, \hat{\mu}_- )^\T \widehat{\blambda} \big\{1 + o_P(1)\big\} \\
&=& 2n \widehat{\blambda}^\T \big\{ \bU_n - \widehat{\bW}_n (\widehat{\mu}_- - \mu_-) \big\} - n \widehat{\blambda}^\T \widehat{\bV}_n \widehat{\blambda} \big\{1 + o_P(1)\big\} \\
&=& nh \bU_n^\T \bV^{-1/2} \bP \bV^{-1/2} \bU_n \big\{1 + o_P(1)\big\}.
\end{eqnarray*}
By  (\ref{proof-3-0}) and the fact that $\bP$ is a rank-1 projection matrix, we conclude that
$$
l_{\SRD,\TR}(\tau_{\SRD}) \stackrel{\mathcal{D}}{\to} \chi^2_1.
$$

In the following, we will demonstrate that $\hat{\mu}_- - \mu_- = O_P\{(nh)^{-1/2}\}$. To do so, we calculate the second derivative of $l_{\SRD,\TR}(\tau_{\SRD}, a)$ with respect to $a$. This derivative is strictly positive definite with probability approaching 1 within a neighborhood set $\mathcal{A}_{\delta} = \{a: |a - \mu_-| \le \delta\} \subseteq \mathcal{A}$, where $\delta$ is a small positive constant. Consequently, $l_{\SRD,\TR}(\tau_{\SRD}, a)$ is a convex function of $a$ over $\mathcal{A}_{\delta}$. Note that $l_{\SRD,\TR}(\tau_{\SRD}, \mu_-) \ge l_{\SRD,\TR}(\tau_{\SRD}, \widehat{\mu}_-)$. It suffices to show that for any constant $\varepsilon > 0$, there exists a sufficiently large constant $C > 0$ such that both $l_{\SRD,\TR}(\tau_{\SRD}, \mu_- - C(nh)^{-1/2})$ and $l_{\SRD,\TR}(\tau_{\SRD}, \mu_- + C(nh)^{-1/2})$ are larger than $l_{\SRD,\TR}(\tau_{\SRD}, \mu_-)$ with high probability, at least $1 - \varepsilon$.

Observe that $l_{\SRD,\TR}(\tau_{\SRD}, \mu_-) = 2 \log \{1 + \widehat{\blambda}_0^\T \bZ_i^{\star}(\tau_{\SRD}, \mu_-) \}$, where $\widehat{\blambda}_0$ satisfies:
$$
\sum_{i=1}^n \frac{\bZ_{i}(\tau_{\SRD}, \mu_-)}{1 + \widehat{\blambda}_0^\T \bZ_i^{\star}(\tau_{\SRD}, \mu_-)} = 0.
$$
It can be shown that $\widehat{\blambda}_0 = \widehat{\bV}_n^{-1} \bU_n + o_P(n^{-1/2}h^{1/2})$ and $\max_{1 \le i \le n} |\widehat{\blambda}_0^\T \bZ_i^{\star}(\tau_{\SRD}, \mu_-)| = o_P(1)$. Hence, we have
$$
l_{\SRD,\TR}(\tau_{\SRD}, \mu_-) = n \bU_n^\T \widehat{\bV}_n^{-1} \bU_n \big\{1 + o_P(1)\big\}.
$$
Analogously, $l_{\SRD,\TR}(\tau_{\SRD}, \mu_- + C(nh)^{-1/2}) = 2 \log \{1 + \widehat{\blambda}_r^\T \bZ_i^{\star}(\tau_{\SRD}, \mu_- + C(nh)^{-1/2}) \}$, where $\widehat{\blambda}_r$ satisfies:
$$
\sum_{i=1}^n \frac{\bZ_{i}(\tau_{\SRD}, \mu_- + C(nh)^{-1/2})}{1 + \widehat{\blambda}_r^\T \bZ_i^{\star}(\tau_{\SRD}, \mu_- + C(nh)^{-1/2})} = 0.
$$
Furthermore, $\widehat{\blambda}_r = \widehat{\bV}_n^{-1} \{\bU_n - C (nh)^{-1/2} \widehat{\bW}_n\} + o_P(n^{-1/2}h^{1/2})$ and $\max_{1 \le i \le n} |\widehat{\blambda}_r^\T \bZ_i^{\star}(\tau_{\SRD}, \mu_- + C(nh)^{-1/2})| = o_P(1)$. Therefore,
\begin{eqnarray*}
&& l_{\SRD,\TR}(\tau_{\SRD}, \mu_- + C(nh)^{-1/2}) \\
&=& n \Big\{ \bU_n - C (nh)^{-1/2} \widehat{\bW}_n \Big\}^\T \widehat{\bV}_n^{-1} \Big\{ \bU_n - C (nh)^{-1/2} \widehat{\bW}_n \Big\} \big\{1 + o_P(1)\big\} \\
&=& n \bU_n^\T \widehat{\bV}_n^{-1} \bU_n + C^2 h^{-1} \widehat{\bW}_n^\T \widehat{\bV}_n^{-1} \widehat{\bW}_n - 2 C n (nh)^{-1/2} \widehat{\bW}_n^\T \widehat{\bV}_n^{-1} \bU_n + o_P(1) \\
&=& l_{\SRD,\TR}(\tau_{\SRD}, \mu_-) + C^2 h^{-1} \widehat{\bW}_n^\T \widehat{\bV}_n^{-1} \widehat{\bW}_n - 2 C (n/h)^{1/2} \widehat{\bW}_n^\T \widehat{\bV}_n^{-1} \bU_n + o_P(1).
\end{eqnarray*}
Notice that $(n/h)^{1/2} \widehat{\bW}_n^\T \widehat{\bV}_n^{-1} \bU_n = O_P(1),$~$
h^{-1} \widehat{\bW}_n^\T \widehat{\bV}_n^{-1} \widehat{\bW}_n = \bW^\T \bV^{-1} \bW \{1 + o_P(1)\},
$~and 
$\bW^\T \bV^{-1} \bW$ is strictly positive definite. By choosing a sufficiently large $C > 0$, we obtain that $l_{\SRD,\TR}(\tau_{\SRD}, \mu_- + C(nh)^{-1/2}) > l_{\SRD,\TR}(\tau_{\SRD}, \mu_-)$ with probability at least $1 - \varepsilon/2$. Similar arguments also yield $l_{\SRD,\TR}(\tau_{\SRD}, \mu_- - C(nh)^{-1/2}) > l_{\SRD,\TR}(\tau_{\SRD}, \mu_-)$ with probability at least $1 - \varepsilon/2$. Hence,
$$
l_{\SRD,\TR}(\tau_{\SRD}, \mu_- \pm C(nh)^{-1/2}) > l_{\SRD,\TR}(\tau_{\SRD}, \mu_-)
$$
with probability at least $1 - \varepsilon$. This implies that $\hat{\mu}_- - \mu_- = O_P\{(nh)^{-1/2}\}$.
The proof for $l_{\SRD,\TR}(\tau_{\SRD})$ is now complete. 

Similar proof techniques can be applied to $l_{\SRD,\DR}(\tau_{\SRD})$ for the sharp RDD, as well as $l_{\FRD,\TR}(\tau_{\FRD})$ and $l_{\FRD,\DR}(\tau_{\FRD})$ for the fuzzy RDD. See also the proof of Theorem 3.1 in \cite{Otsu2015}. We thus complete the proofs of Theorems~\ref{th3} and \ref{th4}.
$\hfill\Box$

\onehalfspacing
\bibliographystyle{dcu}
\bibliography{ref}

\end{document}